\documentclass[a4paper]{article}

\usepackage{amsthm}
\usepackage{amssymb}
\usepackage{latexsym}
\usepackage{longtable}
\usepackage{epsfig}
\usepackage{amsmath}
\usepackage{hhline}
\usepackage{epic}

\newtheorem{thm}{Theorem}[section]
\newtheorem{defn}[thm]{Definition}

\newtheorem{lem}[thm]{Lemma}

\newtheorem{prop}[thm]{Proposition}

\begin{document}

\title{Generalized Blocks of Unipotent Characters in the Finite General Linear Group}
\author{Jean-Baptiste Gramain \\ \'Ecole Polytechnique F\'ed\'erale de Lausanne\\Lausanne, Switzerland \\jean-baptiste.gramain@epfl.ch}
\date{July, 2007}
\maketitle

\begin{abstract}
In a paper of 2003, B. K\"ulshammer, J. B. Olsson and G. R. Robinson
defined $\ell$-blocks for the symmetric groups, where $\ell >1$ is
an arbitrary integer, and proved that they satisfy an analogue of
the Nakayama Conjecture. Inspired by this work and the definitions
of generalized blocks and sections given by the authors, we give in
this paper a definition of $d$-sections in the finite general linear
group, and construct $d$-blocks of unipotent characters, where $d
\geq 1$ is an arbitrary integer. We prove that they satisfy one
direction of an analogue of the Nakayama Conjecture, and, in some
cases, prove the other direction. We also prove that they satisfy an
analogue of Brauer's Second Main Theorem.
\end{abstract}

\section{Introduction, generalized blocks}

Throughout this paper, we let $n$ be a positive integer and $q$ be a
power of a prime $p$. We let $V$ be an $n$-dimensional vector space
over a finite field ${\bf F}_q$ with $q$ elements. We will work in
the finite general linear group $G=GL(V)$. It will sometimes be
convenient to choose a basis for $V$, and then identify $G$ with the
group $GL(n,q)$ of invertible $n$ by $n$ matrices with entries in
${\bf F}_q$. The irreducible complex characters of $G$ have been
described by J. A. Green in \cite{Green}, using deep combinatorial
arguments. Then, using in particular the Deligne-Lusztig theory, P.
Fong and B. Srinivasan have classified the blocks of $G$ (cf
\cite{Fo-Sri}). The unipotent characters of $G$ are parametrized by
the partitions of $n$. It turns out that, if $r$ is a prime not
dividing $q$, then two unipotent characters belong to the same
$r$-block of $G$ if and only if the partitions labeling them have
the same $e$-core, where $e$ is the multiplicative order of $q$
modulo $r$. This result is shown using analogues of the
Murnaghan-Nakayama rule for irreducible characters of $G$. Our aim
is to use these analogues to obtain properties of some
{\emph{generalized blocks}} we define in $G$, and which don't depend
on any prime. We construct \emph{unipotent blocks} for $G$ which
satisfy one direction of an analogue of the Nakayama Conjecture. In
some cases, we also prove that they satisfy both directions.

\bigskip
The concept of {\emph{generalized blocks}} was introduced by B.
K\"{u}lshammer, J. B. Olsson and G. R. Robinson (cf \cite{KOR}).
Take any finite group $H$, and take a union ${\cal C}$ of conjugacy
classes of $H$ containing the identity. We can consider the
restriction to ${\cal C}$ of the ordinary scalar product on
characters of $H$. We denote by Irr$(H)$ the set of complex
irreducible characters of $H$. For $\chi, \, \psi \in
{\mbox{Irr}}(H)$, we let
$$<\chi, \, \psi>_{\cal C}:= \displaystyle \frac{1}{|H|} \sum_{h \in {\cal C}} \chi(h) \psi(h^{-1}).$$
Then $\chi$ and $\psi$ are said to be {\emph{directly}} ${\cal
C}${\emph{-linked}} if $<\chi, \, \psi>_{\cal C} \neq 0$, and
{\emph{orthogonal across}} ${\cal C}$ otherwise. Then direct ${\cal
C}$-linking is a reflexive (since $1 \in {\cal C}$) and symmetric
binary relation on Irr$(H)$. Extending it by transitivity, we obtain
an equivalence relation (called ${\cal C}${\emph{-linking}}) on
Irr$(H)$ whose equivalence classes are called the ${\cal C}$-blocks.
Note that, since they are orthogonal on the whole of $H$, two
{\bf{distinct}} irreducible characters are directly linked across
${\cal C}$ if and only if they are directly linked across $H
\setminus {\cal C}$. Note also that, if we take ${\cal C}$ to be the
set of $r$-regular elements of $H$, for some prime $r$, then the
${\cal C}$-blocks are just the $r$-blocks (cf e.g. \cite{Navarro}).

\medskip
In \cite{KOR}, the authors have defined $\ell$-blocks for the
symmetric groups, where $\ell \geq 2$ is any integer. To obtain
this, they take ${\cal C}$ to be the set of $\ell${\emph{-regular}}
elements, i.e. none of whose cycle has length divisible by $\ell$
(in particular, if $\ell$ is a prime $r$, then the $\ell$-blocks are
just the $r$-blocks). The $\ell$-blocks thus obtained satisfy an
analogue of the Nakayama Conjecture: two irreducible characters
$\varphi_{\lambda}$ and $\varphi_{\mu}$ of the symmetric group $S_n$
(where $\lambda$ and $\mu$ are partitions of $n$) belong to the same
$\ell$-block if and only if $\lambda$ and $\mu$ have the same
$\ell$-core.

Following this work, A. Mar\'oti studied generalized blocks in the
alternating groups, and proved that, if $\ell$ is 2 or any odd
integer greater than 1, then the $\ell$-blocks of the alternating
groups also satisfy an analogue of the Nakayama Conjecture (cf
\cite{Attila}).

\medskip
In the case of symmetric groups, the argument for one of the
directions goes as follows. Any element of a symmetric group can be
written uniquely as a disjoint product of an $\ell$-regular element
and an $\ell$-element (i.e. all of whose cycles have length
divisible by $\ell$). Two elements are said to belong to the same
$\ell${\emph{-section}} if their $\ell$-parts are conjugate. The
$\ell$-section of the identity is thus the set of $\ell$-regular
elements. The $\ell$-sections have properties very similar to those
of ordinary prime sections. Even though we don't need sections in
order to construct generalized blocks, they are important when one
wants to discuss properties of blocks such as analogues of Brauer's
Second Main Theorem.

If $x$ is an $h$-cycle of the symmetric group $S_n$ and $y \in S_n$
fixes the points moved by $x$ (and can thus be seen as an element of
$S_{n-h}$), then the Murnaghan-Nakayama Rule allows us the write the
value of a character $\varphi_{\lambda}$ of $S_n$ at $xy$ as a
linear combination of the $\varphi_{\mu}(y)$'s, where ${\mu}$ runs
through the set of partitions of $n-h$ which can be obtained by
removing an $h$-hook from $\lambda$.

Repeated use of this rule allows us to remove all of the $\ell$-part
of an element by removing a sequence of $\ell$-hooks from the
partition labeling the character we consider. An induction argument
then proves that, if two characters of $S_n$ are linked across any
$\ell$-section, then the partitions labeling them must have the same
$\ell$-core.

\bigskip
Now, in $G=GL(V)$, we don't have such a nice cycle structure
anymore. However, we can use the {\emph{rational canonical form}},
which is based on the {\emph{elementary divisors}} of a given
element $g$. These are polynomials over ${\bf F}_q$, and they
correspond to the cyclic subspaces of $V$ under the action of $g$.
Then, very much like in the symmetric group, the Murnaghan-Nakayama
Rule for unipotent characters allows us to relate the "removal" of a
cyclic subspace associated to a polynomial of degree $d$ to the
removal of $d$-hooks from the partition labeling the unipotent
character we consider.

The definitions of $d${\emph{-element}} and $d${\emph{-regular
element}} we give are such that we can apply the Murnaghan-Nakayama
Rule and show by an induction argument that, if two unipotent
characters are linked across any $d${\emph{-section}}, then the
partitions labeling them have the same $d$-core.

\bigskip
In section 2, we present the conjugacy classes of $GL(V)$ and the
{\emph{primary decomposition}}, which is the basis for the
definitions we introduce. In section 3, we define the $d$-sections
we work on, and the generalized blocks they define. Section 4 is
devoted to the analogue of the Nakayama Conjecture for unipotent
$d$-blocks. We prove that one direction is true in general, and that
the other is true for $d=1$ and, for any $d>0$, for blocks of weight
at most $(d+1)/2$. We conjecture that, in fact, these restrictions
are not necessary. Finally, in section 5, we discuss an analogue of
Brauer's Second Main Theorem for unipotent $d$-blocks.

\section{Conjugacy classes, primary decomposition}

\subsection{Rational canonical form}
We first introduce the theory of elementary divisors and the
rational canonical form in $G$, which gives a parametrization of the
conjugacy classes. For the results we give in this section, we refer
to Green \cite{Green}.

For any non-negative integer $k$, we will write $\nu \vdash k$ to
say that $\nu$ is a partition of $k$, and we will write $|\nu|=k$.
We denote by ${\cal P}$ the set of all partitions. We write ${\cal
F}=\{f_i , \, i \in {\cal I} \}$ the set of irreducible, monic
polynomials distinct from $X$ and of degree at most $n$ in ${\bf
F}_q[X]$ (so that ${\cal I}$ is finite). For $f \in {\cal F}$, we
write $\delta (f)$ for the degree of $f$. Then the conjugacy classes
of $G$ are parametrized by the set of partition-valued functions $
\zeta \colon {\cal F} \longrightarrow {\cal P}$ such that $\sum_{f
\in {\cal F} } | \zeta (f) | \delta (f) =n$. For such a $\zeta$, the
corresponding conjugacy class of $G$ will be written $(f_1^{\nu_1},
\, \ldots , \, f_r^{\nu_r})$, where $\{f_1, \, \ldots, \, f_r \}=\{
f \in {\cal F} , \, | \zeta (f)| \neq 0 \}$, and, for $1 \leq i \leq
r$, $\nu_i = \zeta (f_i)$. We will sometimes also use the notation
$(f^{\zeta (f)})_{f \in {\cal F}}$ (thus allowing $|{\zeta(f)}|$ to
be 0).

Take a conjugacy class $c=(f_1^{\nu_1}, \, \ldots , \, f_r^{\nu_r})$
of $G$, and $g \in c$. Then, the characteristic polynomial of $g$
over ${\bf F}_q$ is
$$Char(g)=\prod_{i=1}^{r} f_i^{| \nu_i |}.$$
Writing, for each $i$, $\nu_i=(\lambda_{i,1} \geq \lambda_{i,2} \geq
\ldots \geq \lambda_{i,s} )$ (where $s$ can be chosen to be big
enough for all $f_i's$ appearing in $Char(g)$ (and even for all $g
\in G$), by taking for example $s=n$), the minimal polynomial of $g$
over ${\bf F}_q$ is
$$Min(g)=\prod_{i=1}^{r} f_i^{\lambda_{i,1}}.$$
For $k=1, \, \ldots , s$, the polynomial $ED_k(g)=\prod_{i=1}^r
f_i^{\sum_{j=1}^{k} \lambda_{i,j}}$ is the $k$-th elementary divisor
of $g$ over ${\bf F}_q$.

For any irreducible monic polynomial $f$ over ${\bf F}_q$, we let
$U(f)$ be the companion matrix of $f$, i.e. if
$f(X)=X^d-a_{d-1}X^{d-1}- \cdots - a_0$, then
$$U(f)=U_1(f)= \left( \begin{array}{ccccc} 0 & 1 & & & \\ \vdots & \ddots
& \ddots &
(0) & \\ \vdots & & \ddots & \ddots & \\ 0 & \cdots & \cdots & 0 & 1 \\
a_0 & a_1 & \cdots & \cdots & a_{d-1} \end{array} \right).$$ For any
positive integer $\lambda$, we write
$$U_{\lambda}(f)= \left( \begin{array}{ccccc} U(f) & I_d & & & \\  & U(f)
& I_d &  &
\\  & & \ddots & \ddots & \\  &  &  & \ddots & I_d \\  &  &  &  & U(f)
\end{array}
\right),$$ where $I_d$ is the $d$ by $d$ identity matrix, and there
are $\lambda$ blocks $U(f)$ on the diagonal. Note that
$U_{\lambda}(f)$ is equivalent over ${\bf F}_q$ to the companion
matrix of $f^{\lambda}$. Finally, for $\nu = (\lambda_1 \geq  \ldots
\geq \lambda_s)$ any partition (of an integer $k$ say), we write
$$U_{\nu}(f)= \left( \begin{array}{ccc} U_{\lambda_1}(f) & & \\  & \ddots
& \\  & & U_{\lambda_s}(f) \end{array} \right),$$where $U_0(f)$ has
to be seen as the empty matrix. If $g \in c=(f_1^{\nu_1}, \, \ldots
, \, f_r^{\nu_r})$, then, in {\bf{any}} matrix representation of
$G$, $g$ is equivalent over ${\bf F}_q$ to
$$U_{(f_1^{\nu_1}, \, \ldots , \, f_r^{\nu_r})}= \left(  \begin{array}{ccc}
U_{\nu_1}(f_1) & & \\
& \ddots & \\  & & U_{\nu_r}(f_r) \end{array} \right).$$ This is the
rational canonical form of $g$ over ${\bf F}_q$.

\subsection{The Jordan decomposition}
An element $g \in G$ is {\emph{semisimple}} if and only if it is
diagonalizable over an algebraic extension of ${\bf F}_q$ (i.e. some
${\bf F}_{q^k}$, $k \geq 1$), i.e. $Min(g)$ splits over ${\bf
F}_{q^k}$ and has only simple roots. Equivalently, $g \in GL(V)$ is
semisimple if and only if there is a basis of $V$ consisting of
eigenvectors of $g$. Let $c$ be the conjugacy class of $g$. If
$c=(f_1^{\nu_1}, \, \ldots , \, f_r^{\nu_r})$, then $g$ is
semisimple if and only if, for all $1 \leq i \leq r$, we have
$\nu_i=(1, \, \ldots ,1)$.

Choose some basis for $V$. We write $\Delta ((q_i))$ for the
diagonal matrix with diagonal blocks the $q_i$'s. Suppose
$g=\Delta(( U_{\nu_i}(f_i) ) )$. Then $g=g_S+g_N=\Delta((
U_{\tilde{\nu}_i}(f_i) ) )+\Delta(( \tilde{I}_{\nu_i}(f_i) ) )$,
where $\tilde{\nu_i}=(1, \, \ldots ,1) \vdash | \nu_i |=k_i$ and
$$\tilde{I}_{\nu_i}(f_i)= \left( \begin{array}{ccccc} (0) & I_{d_i} & & &
\\  & (0)
& I_{d_i} &  & \\  & & \ddots & \ddots & \\  &  &  & \ddots & I_{d_i} \\
&  &  &  & (0) \end{array} \right),$$ with $| \nu_i |=k_i$ diagonal
$d_i$ by $d_i$ blocks, where $f_i$ has degree $d_i$.

Thus $g_S=\Delta (( U_{k_i}(f_i) ) )$ is the semisimple part of $g$,
and the unipotent part of $g$ is $g_U=I_n+g_S^{-1}g_N$.

\medskip
Now, if $\tilde{g}=hgh^{-1}$ for some $h \in G$, then
$\tilde{g}=(hg_Sh^{-1})(hg_Uh^{-1})$, and we see that the semisimple
and unipotent parts of $\tilde{g}$ are $hg_Sh^{-1}$ and $hg_Uh^{-1}$
respectively.

\subsection{Primary decomposition}
We first mention a fact about centralizers in $GL(V)$. Take $g \in
GL(V)$. If, in some matrix representation corresponding to the
decomposition $V=V_1 \oplus V_2$, we have
$$g=\left( \begin{array}{cc} g_1 & \\ & g_2 \end{array} \right)$$
with, for $i \in \{ 1,2 \}$, $g_i \in G_i=GL(V_i)$, and
gcd$(Min(g_1),Min(g_2))=1$, then
$$C_G(g)  = \left\{ \left( \begin{array}{cc} a & \\ & b \end{array} \right), \; a \in C_{G_1}(g_1), \, b \in C_{G_2}(g_2) \right\},$$
and thus $C_G(g)  \cong C_{G_1}(g_1) \times C_{G_2}(g_2)$ via $h
\longmapsto (h|_{V_1}, \, h|_{V_2}) $ (which is independant on the
matrix representation).

\medskip
We also give the following lemma
\begin{lem}
Let $g \in G=GL(V)$. Suppose $V_1$ and $V_2$ are $g$-stable
subspaces of $V$ such that $V=V_1 \oplus V_2$. If $g|_{V_1} \in
(f^{\zeta_1(f)})_{f \in {\cal F}}$ and  $g|_{V_2} \in (
f^{\zeta_2(f)})_{f \in {\cal F}}$ for some $\zeta_1, \, \zeta_2
\colon {\cal F} \longrightarrow {\cal P}$, then $g \in
(f^{\zeta_1(f) \cup \zeta_2(f)} )_{f \in {\cal F}}$, where, by $
\cup $, we denote the concatenation of partitions (i.e., for
partitions $\mu$ and $\nu$, the components of the partition $\mu
\cup \nu$ are those of $\mu$ together with those of $\nu$. If $\mu
\vdash k$ and $\nu \vdash \ell$, then $\mu \cup \nu \vdash k+\ell$.)
\end{lem}

\begin{proof}
Taking any bases for $V_1$ and $V_2$ (which then add to give a basis
of $V$), we can write $$g = \left( \begin{array}{cc} g_1 & \\ & g_2
\end{array} \right)$$ (where $g_1$ (resp. $g_2$) is the matrix of
$g|_{V_1}$ (resp. $g|_{V_2}$)). The idea is that we can obtain the
rational canonical form of $g$ by reducing to this form $g_1$ and
$g_2$. There exist $h_1 \in GL(V_1)$ and $h_2 \in GL(V_2)$ such that
$$g=\left( \begin{array}{cc} h_1 & \\ & h_2 \end{array} \right)^{-1}
\left( \begin{array}{cc} U_{(f^{\zeta_1(f)})_{f \in {\cal F}}} & \\
& U_{( f^{\zeta_2(f)})_{f \in {\cal F}}}
\end{array} \right) \left( \begin{array}{cc} h_1 & \\ & h_2
\end{array} \right).$$ Thus, for some permutation matrix $P$,
$$g=\left( \begin{array}{cc} h_1 & \\ & h_2 \end{array} \right)^{-1} P^{-1}
U_{(f^{\zeta_1(f) \cup \zeta_2(f)} )_{f \in {\cal F}}} P \left(
\begin{array}{cc} h_1 & \\ & h_2
\end{array} \right),$$ so that $g \in (f^{\zeta_1(f) \cup \zeta_2(f)} )_{f \in {\cal F}}$.
\end{proof}

We now turn to the primary decomposition of elements of $G$. Take
\\
$g \in GL(V)$ and suppose $g \in (f_1^{\nu_1}, \,  \ldots , \,
f_r^{\nu_r})$. Then there exists a unique decomposition $V = V_1
\oplus \cdots \oplus V_r$, where the $V_i$'s are $g$-stable
subspaces of $V$ and, for $1 \leq i \leq r$, $g|_{V_i} \in
(f_i^{\nu_i}) \subset GL(V_i)$. For each $1 \leq i \leq r$, $V_i$ is
given by
$$V_i= \{ v \in V \, | \, f_i^k(g)v=0 \; \mbox{for some} \; k >0 \}.$$
We have $C_G(g) \cong C_{GL(V_1)}(g|_{V_1}) \times \cdots \times
C_{GL(V_r)}(g|_{V_r})$. Then there exists a unique writing $g=g_1
\ldots g_r$, where, for each $1 \leq i \leq r$, $g_i \in GL(V)$,
$V_i$ is $g_i$-stable and $g_i|_{V_j}=1$ for all $j \neq i$. Indeed,
we must have $g_i|_{V_i}=g|_{V_i}$ for each $1 \leq i \leq r$, and
$g_i$ is uniquely determined by this and its properties listed
before. Furthermore, the $g_i$'s are pairwise commuting elements. We
say that $g_1 \ldots g_r$ is the \emph{primary decomposition} of g.
The $g_i$'s are the \emph{primary components} of $g$.

More generally, an element of $G$ is said to be \emph{primary} if
its characteristic polynomial is divisible by at most one
irreducible polynomial distinct from $X-1$. We have the following
(cf Fong-Srinivasan \cite{Fo-Sri}, Proposition (1A))
\begin{prop}
Suppose that $h$ is a semisimple primary element of some general
linear group $GL(m,q)$, and that $h \in (f^{\nu})$ for some $f \neq
X-1$. Writing $d$ the degree of $f$, we thus have $m=kd$, where
$\nu=(1, \ldots,1) \vdash k$. Then
$$C_{GL(kd,q)}(h) \cong GL(k,q^d).$$
\end{prop}

\medskip
From the primary decomposition of $g \in G$, we deduce the
following: if ${\cal F}$ is the disjoint union of ${\cal F}_1$ and
${\cal F}_2$, then there is a unique decomposition $V=V_1 \oplus
V_2$ where $V_1$ and $V_2$ are $g$-stable, $g|_{V_1} \in
(f^{\nu(f)})_{f \in {\cal F}_1}$ and $g|_{V_2} \in (f^{\nu(f)})_{f
\in {\cal F}_2}$. Then $g$ has a unique decomposition
$g=g_1g_2=g_2g_1$ where $g_1,g_2 \in GL(V)$, $g_1|_{V_2}=1$ and
$g_2|_{V_1}=1$. Indeed, under these hypotheses, necessarily, for $i
\in \{ 1,2 \}$, $V_i$ is $g_i$-stable and $g_i|_{V_i}=g|_{V_i}$. We
then have $C_G(g) \cong C_{GL(V_1)}(g|_{V_1}) \times
C_{GL(V_2)}(g|_{V_2})$.

\section{Sections, blocks}

\subsection{Sections}

The idea of the following definitions is, in the rational canonical
form we gave above, to isolate blocks corresponding to irreducible
polynomials whose degree is divisible by a given (positive) integer
$d$. We let, writing $\delta (f)$ the degree of any polynomial $f$,
$${\cal F}_{d}=\{ f \in {\cal F} \, | \, f \neq X-1 \; \mbox{and} \; d | \delta (f) \} \; \mbox{and}
\; {\cal F}_{d}^0={\cal F} \setminus {\cal F}_{d}.$$

We define the following union of conjugacy classes of $G$:
$${\cal X}_d= \{ g \in (f^{\nu(f)})_{f \in {\cal F}_d \cup \{ X-1 \} } \, | \, \nu(X-1)=(0) \, \mbox{or} \,(1, \ldots, 1)
\}.$$ Note that we have $1 \in {\cal X}_d$.

\smallskip
For $x \in {\cal X}_d$, we let ${\cal F}_x={\cal F}_d$ and ${\cal
F}_x^0={\cal F}_d^0$. Then, for each $x \in {\cal X}_d$, the set
${\cal F}$ is the disjoint union of ${\cal F}_x$ and ${\cal F}_x^0$.
For each $x \in {\cal X}_d$, there exists a unique decomposition
$V=V_x \oplus V_x^0$ such that $V_x$ is $x$-stable, $x |_{V_x^0}=1$
and $x|_{V_x} \in (f^{\nu(f)})_{f \in {\cal F}_x}$. We then have
$C_G(x) \cong C_{GL(V_x)}(x|_{V_x}) \times GL(V_x^0)$, and we write
$C_G(x)=H_x \times H_x^0$. We let
$${\cal Y}_d(x)= \{ y \in GL(V) \, | \, V_x^0 \, \mbox{is} \, y \, \mbox{-stable}, \, y|_{V_x}=1 \,
\mbox{and} \, y|_{V_x^0} \in ( f^{\nu(f)} )_{f \in {\cal F}_x^0}
\}.$$ From the definitions, we see that ${\cal Y}_d(x) \subset
C_G(x)$ for any $x \in {\cal X}_d$. We also see, using the remarks
we made on the primary decomposition, that, for any $g \in G$, there
exist {\bf{unique}} $x \in {\cal X}_d$ and $y \in {\cal Y}_d(x)$
such that $g=xy$. Indeed, if $g \in (f_1^{\nu_1}, \ldots,
f_r^{\nu_r})$ and if we write, as above, $V=V_1 \oplus \cdots \oplus
V_r$ and $g =g_1 \ldots g_r$, then we necessarily have $V_x=
\oplus_{i \in I} V_i$, where $I \subset \{ 1, \, \ldots , \, r \}$
is the set of indices $i$ such that $f_i^{\nu_i}$ has the property
defining ${\cal X}_d$, $V_x^0=\oplus_{i \not \in I} V_i$,
$x|_{V_x}=g|_{V_x}$ and $y|_{V_x^0}=g|_{V_x^0}$ (and $x \in {\cal
X}_d$ and $y \in {\cal Y}_d(x)$ are uniquely determined by these
conditions). We have $x=\prod_{i \in I}g_i$ and $y = \prod_{i \not
\in I} g_i$.

These definitions allow us to define $({\cal X}_d,{\cal
Y}_d)$-sections as introduced in section 2 of \cite{KOR}. The
following proposition shows that ${\cal X}_d$ and ${\cal Y}_d(x)$
($x \in {\cal X}_d$) behave in many respects like the set of
$r$-elements ($r$ a prime) and the sets of $r$-regular elements of
their centralizers respectively.

\begin{prop}
Take any integer $d>0$. Then, for any $x \in {\cal X}_d$,

\noindent (i) ${\cal Y}_d(x)$ is a union of conjugacy classes of
$C_G(x)$.

\noindent (ii) For all $y \in {\cal Y}_d(x)$, $C_G(xy) \leq C_G(x)$.

\noindent (iii) For all $g \in G$, ${\cal Y}_d(x^g)={\cal
Y}_d(x)^g$.

\noindent (iv) Two elements of $x {\cal Y}_d(x)$ are $G$-conjugate
if and only if they are $C_G(x)$-conjugate.

\noindent (v) $G=\coprod_{x \in {\cal X}_d/G} \{ (xy)^G, \, y \in
{\cal Y}_d(x)/C_G(x) \}$.
\end{prop}

\begin{proof}

\noindent (i) We have $C_G(x)=H_x \times H_x^0$, and we see $C_G(x)$
as a subgroup of $GL(V_x) \times GL(V_x^0)$. For $y \in {\cal
Y}_d(x)$, we have $y=(y|_{V_x},y|_{V_x^0})=(1,y_x^0)$. Then, for all
$h=(h_x,h_x^0) \in C_G(x)$, $h^{-1}yh=(1,(h_x^0)^{-1}y_x^0h_x^0) \in
{\cal Y}_d(x)$.

\noindent (ii) For all $y \in {\cal Y}_d(x)$, we have $C_G(xy) \cong
C_{GL(V_x)}(x|_{V_x}) \times C_{GL(V_x^0)}(y|_{V_x^0}) \leq
C_{GL(V_x)}(x|_{V_x}) \times GL(V_x^0) \cong C_G(x)$, and $C_G(xy)
\leq C_G(x)$ since the isomorphism on the left is the restriction to
$C_{GL(V_x)}(x|_{V_x}) \times C_{GL(V_x^0)}(y|_{V_x^0})$ of the
isomorphism on the right.

\noindent (iii) Take $g \in G$ and $y \in {\cal Y}_d(x)$. Then $x^g
\in {\cal X}_d$ and ${\cal F}_{x^g}={\cal F}_x$ (and thus ${\cal
F}_{x^g}^0={\cal F}_x^0$). We have $V=g^{-1}V=g^{-1}V_x \oplus
g^{-1} V_x^0$. Furthermore, $x^g|_{g^{-1}V_x^0}=1$, $g^{-1}V_x$ is
$x^g$-stable, and $x^g|_{g^{-1}V_x} \in (f^{\nu})_{f \in {\cal
F}_x}$. Thus $g^{-1}V_x=V_{x^g}$ and $g^{-1}V_x^0=V_{x^g}^0$. Now,
since $y \in {\cal Y}_d(x)$, we see that $y^g|_{g^{-1}V_x}=1$,
$g^{-1}V_x^0$ is $y^g$-stable, and $y^g|_{g^{-1}V_x^0} \in
(f^{\nu})_{f \in {\cal F}_x^0}$. Hence $y^g \in {\cal Y}(x^g)$, and
${\cal Y}_d(x)^g \subset {\cal Y}_d(x^g)$ for all $g\in G$. Then,
for any $g \in G$, ${\cal Y}_d(x^g)^{g^{-1}} \subset {\cal Y}_d(x)$
so that ${\cal Y}_d(x^g) \subset {\cal Y}_d(x)^g$. Hence the result.

\noindent (iv) Suppose that, for some $y,z \in {\cal Y}_d(x)$, there
exists $h \in G$ such that $xy=h^{-1}xzh$. Writing $g=xy$, we also
have $g=x^hz^h$, and $x^h \in {\cal X}_d$ (since ${\cal X}_d$ is a
union of $G$-conjugacy classes) and $z^h \in {\cal Y}_d(x^h)$ (by
(iii)). By the uniqueness of such a writing for $g$, we have $x^h=x$
and $z^h=y$. Hence $h \in C_G(x)$. In particular, $xy$ and $xz$ are
$C_G(x)$-conjugate.

\noindent (v) For any $g \in G$, there exists a unique $x \in {\cal
X}$ such that $g \in x{\cal Y}_d(x)$. Thus $$\begin{array}{cl} G & =
\displaystyle \coprod_{x \in {\cal X}_d} x {\cal Y}_d(x) \\ & =
\displaystyle \coprod_{x \in {\cal X}_d/G} (x {\cal Y}_d(x))^G \;
\mbox{(because of (iii))} \\ & =\displaystyle \coprod_{x \in {\cal
X}/G} \; \;  \bigcup_{y \in {\cal Y}_d(x)} (xy)^G  \\ & =
\displaystyle  \coprod_{x \in {\cal X}_d/G} \; \; \coprod_{y \in
{\cal Y}_d(x)/C_G(x)} (xy)^G, \end{array}$$ this last equality being
true by (iv) (if $(xy)^g=(xz)^{g'}$, then $xy=xz^{g'g^{-1}}$ so
that, by (iv), $y$ and $z$ are $C_G(x)$-conjugate).

\end{proof}

\medskip
Because of the analogy we mentioned before the proposition, we
introduce some terminology. We will refer to the elements of ${ \cal
X}_d$ as $d$-{\emph{elements}}. One set of particular importance is
${\cal Y}_d(1)$. This is the set of elements of $G$ whose minimal
polynomial has no irreducible factor of degree divisible by $d$
(except possibly $X-1$). We call the elements of ${\cal Y}_d(1)$ the
$d$-{\emph{regular elements}} of $G$. The elements of $G \setminus
{\cal Y}_d(1)$ are called $d$-{\emph{singular}}.

\medskip
Note that, with these definitions, a $1$-regular element of $G$ is
an element whose minimal polynomial has no irreducible factor of
degree divisible by 1, except $X-1$. This means that its minimal
polynomial is a power of $X-1$. Hence the 1-regular elements of
$GL(n,q)$ are precisely the unipotent elements.

\medskip
For any $x \in {\cal X}_d$, we call the union of the $G$-conjugacy
classes meeting $x{\cal Y}_d(x)$ the ${\cal Y}_d${\emph{-section}}
of $x$. The above proposition implies that, for each $x \in {\cal
X}_d$, induction of complex class functions gives an isometry from
the space of class functions of $C_G(x)$ vanishing outside $x{\cal
Y}_d(x)$ onto the space of class functions of $G$ vanishing outside
the ${\cal Y}_d$-section of $x$ (cf \cite{KOR}).

\medskip
We remark that the ${\cal Y}_d$-sections of $G$ are quite different
from ordinary prime sections. Take any $1 \neq x \in {\cal X}_d$.
Then, by definition, ${\cal Y}_d(x) \subset {\cal Y}_d(1) \cap
C_G(x)$. However, if $d \neq 1$, then there exists $\lambda \in {\bf
F}_q^{\times}$ such that $\lambda I_n \in {\cal Y}_d(1) \cap C_G(x)$
but $ \lambda I_n \not \in {\cal Y}_d(x)$, so that ${\cal Y}_d(x)
\neq {\cal Y}_d(1) \cap C_G(x)$ (while this equality holds when we
define $({\cal X}_d,{\cal Y}_d)$-sections to be the ordinary
$r$-sections for some prime $r$).

Furthermore, still supposing $d \neq 1$, if $x$ is a (non trivial)
$r$-element of $G$ for some prime $r$, then, most of the time (that
is, when $q-1$ is not a power of $r$), there exists an $r$-regular
element $\lambda \in {\bf F}_q^{\times}$ such that $\lambda I_n \not
\in {\cal Y}_d(x)$, so that $x \lambda I_n \not \in (x {\cal
Y}(x))^G$. But $x \lambda I_n$ belongs to the $r$-section of $x$.
Hence the ${\cal Y}_d$-section of $x$ is {\bf not} a union of
$r$-sections.

\subsection{Blocks}

With the definition of $d$-regular element we have given, we can now
define generalized blocks of irreducible characters, as described in
the introduction, by taking ${\cal C}$ to be the set ${\cal Y}_d(1)$
of $d$-regular elements. For any positive integer $d$, the
$d${\emph{-blocks}} of $G$ are defined to be minimal subsets of
Irr$(G)$ which are orthogonal across $d$-regular elements. Recall
that two characters $\chi, \, \psi \in $Irr$(G)$ are orthogonal
across $d$-regular elements if
$$<\chi, \, \psi>_{{\cal Y}_d(1)}:=\frac{1}{|G|} \sum_{y \in {\cal Y}_d(1)} \chi(y) \overline{\psi(y)} = 0.$$
Otherwise, $\chi$ and $\psi$ are said to be directly linked across
$d$-regular elements.

% If we take a prime $r$, then the $r$-blocks of
%$G$ are minimal subject to separating $r$-regular elements from
%$r$-singular elements (cf e.g.\cite{Navarro}(......)). Equivalently,
%irreducible characters in distinct $r$-blocks are
%{\emph{orthogonal}} across $r$-regular elements, and the $r$-blocks
%are minimal subsets of Irr$(G)$ with respect to this property.

%Following \cite{KOR}, we define {\emph{generalized blocks}} for $G$ by
%orthogonality across $d$-regular elements, where $d$ is now an
%arbitrary integer.

 %Take any integer $d>0$. We define on Irr$(G)$ the relation
%$\sim$ of direct ${\cal Y}_d(1)$-linking: for $\chi, \, \psi \in
%$Irr$(G)$, $\chi \sim \psi$ if and only if
%$$<\chi, \, \psi>_{{\cal Y}_d(1)}:=\frac{1}{|G|} \sum_{y \in {\cal Y}_d(1)} \chi(y) \overline{\psi(y)} \neq 0.$$
%We then say that $\chi$ and $\psi$ are {\emph{directly linked
%across}} $d$-regular elements. Otherwise, $\chi$ and $\psi$ are said
%to be {\emph{orthogonal across}} $d$-regular elements.

%Extending $\sim$ by transitivity, we obtain an equivalence relation
%$\approx$ on Irr$(G)$. We call the equivalence classes of $\approx$
%the $d${\emph{-blocks}} of $G$.

We will particularly consider the restriction of the relation of
${\cal Y}_d(1)$-linking to the subset of unipotent characters. The
unipotent characters of $G$ are the irreducible components of the
permutation representation of $G$ on the cosets of a Borel subgroup
(i.e. the normalizer in $G$ of a Sylow $p$-subgroup of $G$, where
$p$ is the defining characteristic). They are labeled by the
partitions of $n$. For $\lambda \vdash n$, we will write $\chi
_{\lambda}$ the unipotent (irreducible) character of $G$ labeled by
$\lambda$. We write Unip$(G)=\{ \chi_{\lambda}, \, \lambda \vdash n
\}$. The equivalence classes of Unip$(G)$ modulo the restriction of
${\cal Y}_d(1)$-linking will be called {\emph{unipotent}}
$d${\emph{-blocks}} of $G$. It is clear that, for any $d$-block $B$
of $G$, $B \cap \mbox{Unip}(G)$ is a {\bf{union}} of unipotent
$d$-blocks.

\section{Nakayama Conjecture for unipotent blocks}

\subsection{Murnaghan-Nakayama Rule for unipotent characters}

Pick $g \in G$, and write $Char(g)= \prod_i f_i^{k_i}$ and the
corresponding decomposition $g = \prod_i g_i$. Then pick $i_0$, and
write $g=\rho \sigma$, where
$$\rho=g_{i_0} \; \mbox{and} \; \sigma=\prod_{i \neq i_0}g_i.$$
Now fix some matrix representation of $G$, and see $G$ as $GL(n,q)$.
Writing $d$ the degree of $f_{i_0}$, $m=k_{i_0}d$, and $l=n-m$, we
have, writing $\sim$ for equivalence of matrices over ${\bf F}_q$,
$$g \sim \left( \begin{array}{ccc} \ddots & & \\ & U_{\nu_i}(f_i) & \\ & & \ddots \end{array} \right),$$
$$\rho \sim \left( \begin{array}{cc} I_l &  \\  & U_{\nu_{i_0}}(f_{i_0})  \end{array} \right)$$
(and we may consider that $\rho \in G_0=GL(m,q)$), and
$$\sigma \sim \left( \begin{array}{cccc}  () & & & \\ & \ddots & & \\ & & () & \\ & & & I_m \end{array} \right)$$
(and we may consider that $\sigma \in G_1=GL(l,q)$).

Then, using the results on the Jordan decomposition given in the
first section, we see that the semi-simple part $\rho_S$ of $\rho$
is equivalent to
$$\left( \begin{array}{cc} I_l & \\ & U_{k_{i_0}}(f_{i_0}) \end{array} \right).$$
We have $\rho_S \in GL(n,q)$ and $Char( \rho_S)=f_{i_0}^{k_{i_0}}
(X-1)^l$, and we may consider that $\rho_S \in G_0=GL(m,q)$ and
$Char( \rho_S)=f_{i_0}^{k_{i_0}}$.

We have $C_G(\rho_S)=H=H_0 \times H_1$, where $H_1 \cong GL(l,q)$
and \\
$H_0=C_{GL(m,q)}(\rho_S) \cong GL(k_{i_0},q^d) $.

\medskip
We find in \cite{Fo-Sri} an analogue of the Murnaghan-Nakayama Rule
which applies to unipotent characters of $G$. More precisely, if
$\chi_{\nu} \in$ Unip$(G)$, then the following theorem gives us
information on the values of a class function $\chi^{\nu}$ which is,
{\bf{up to a sign}}, the same as $\chi_{\nu}$ (it is this class
function, rather than the unipotent character, which appears
naturally in the Deligne-Lusztig theory).

\begin{thm}[Murnaghan-Nakayama Rule] (\cite{Fo-Sri}, Theorem (2G))

Let $g \in G$, and $\rho$ and $\sigma$ be as above, and let $\nu
\vdash n$. Then
$$\chi^{\nu}(\rho \sigma)=\sum_{\lambda \in {\cal L}_{\nu}} a_{\nu \lambda}^{\rho} \chi^{\lambda}(\sigma),$$
where ${\cal L}_{\nu}$ is the set of partitions $\lambda$ of $l$
which can be obtained from $\nu$ by removing $k_{i_0}$ $d$-hooks,
and $a_{\nu \lambda}^{\rho} \in {\bf Z}[q^d]$.

\noindent If ${\cal L}_{\nu} = \emptyset$, then $\chi^{\nu}(\rho
\sigma)=0$.

\noindent We have $a_{\nu \lambda}^{\rho} \neq 0$ for $\lambda \in
{\cal L}_{\nu}$.

\smallskip
\noindent (The coefficients of $a_{\nu \lambda}^{\rho}$ depend on
the characters of the symmetric group $S_{k_{i_0}}$ and the Green
functions of $GL(k_{i_0},q^d) \cong H_0$ (applied to the unipotent
part $\rho_U$ of $\rho$), and all the non-zero coefficients of
$a_{\nu \lambda}^{\rho}$ have the same sign).

\end{thm}

{\bf Remark:} it is easy to see that, if $a_{\nu \lambda}^{\rho}
\neq 0$, then $\nu$ and $\lambda$ have the same $d$-core.

\bigskip
The idea is to use this theorem recursively so as to be able to
obtain information about $<\chi^{\lambda}, \, \chi^{\mu}>_{x{\cal
    Y}_d(x)}$, for $\lambda, \, \mu \vdash n$ and $x \in {\cal X}_d$.
     We first use it to obtain a formula for $\chi^{\mu}(xy)$, where $\mu
\vdash n$ and  $y \in {\cal Y}_d(x)$. We take $1 \neq x \in {\cal
X}_d$. Suppose $x \in c_x=( f_1^{\lambda_1}, \, \ldots , \,
f_s^{\lambda_s} )$ where, for each $1 \leq i \leq s$, $\lambda_i
\vdash k_i$, and $d | \delta (f_i)$ or $f_i=X-1$. If $f_i \neq X-1$,
we let $\delta (f_i)=m_id$. From the definition of ${\cal X}_d$, we
see that, in the primary decomposition of $x$, we may omit the
component corresponding to $X-1$, because it is necessarily the
identity. We relabel $f_1, \, \ldots , \, f_r$ the $f_i$'s which are
distinct from $X-1$, and we write $x=x_1 \ldots x_r$, where each
$x_i$ has exactly one elementary divisor distinct from $X-1$
(namely, $f_i$), and with the same multiplicity as in $Char(x)$. We
will say that $x$ has $d${\emph{-type}} $(k_1m_1, \ldots ,k_rm_r)$.

By repeated use of the Murnaghan-Nakayama Rule, we obtain that, for
any $y \in {\cal Y}_d(x)$ and for $\mu \vdash n$, we have
$$\begin{array}{cl}\chi^{\mu}(xy) & = \sum_{\mu_1 \in {\cal L}_{\mu}} a_{\mu \mu_1}^{x_1}
\chi^{\mu_1} (x_2 \ldots x_ry) \\ &  = \sum_{\mu_1 \in {\cal L}_{\mu}} a_{\mu \mu_1}^{x_1}
 \left( \sum_{\mu_2 \in {\cal L}_{\mu_1}} a_{\mu_1 \mu_2}^{x_2} \chi^{\mu_2} (x_3 \ldots x_ry)
  \right) \\ &  = \sum_{\mu_1 \in {\cal L}_{\mu}} \ldots   \sum_{\mu_r \in {\cal L}_{\mu_{r-1}}}
   a_{\mu \mu_1}^{x_1} \ldots  a_{\mu_{r-1} \mu_r}^{x_r} \chi^{\mu_r} (y) \end{array},$$
which can be written
$$\chi^{\mu}(xy)=\displaystyle \sum_{\lambda \in {\cal L}_{(k_1m_1, \ldots,k_rm_r)d}^{\mu}}
\alpha_{\mu \lambda}^{(x_1,\ldots,x_r)} \chi^{\lambda} (y),$$ where
the $\alpha_{\mu \lambda}^{(x_1,\ldots,x_r)}$'s are integers and
${\cal L}_{(k_1m_1, \ldots,k_rm_r)d}^{\mu}$  is the set of
partitions of $n-(\sum_i k_im_i)d$ which can be obtained from $\mu$
by removing $k_1$ $(m_1d)$-hooks, then $k_2$ $(m_2d)$-hooks, \ldots
, and finally $k_r$ $(m_rd)$-hooks. We will call such a sequence of
removals a $(k_1m_1, \ldots , k_rm_r)d${\emph{-path}} from $\mu$ to
$\lambda$.

Note that, in this sum, each $\lambda$ can appear several times, if
there is more than one $(k_1m_1, \ldots, k_rm_r)d$-path from $\mu$
to $\lambda$. Note also that, in the right side of this equality,
$y$ has implicitely been seen as an element of $GL(l,q)$, where
$l=n-(\sum_i k_im_i)d$.

\smallskip
\noindent If ${\cal L}_{(k_1m_1, \ldots,k_rm_r)d}^{\mu} =
\emptyset$, then $\chi^{\mu}(xy)=0$.

\noindent
We have $ \alpha_{\mu \lambda}^{(x_1,\ldots,x_r)} \in {\bf
Z}[q^d]$, and, if we separate the possibly multiple occurences of
each $\lambda$ in the sum, then, for $\lambda \in {\cal L}_{(k_1m_1,
\ldots,k_rm_r)d}^{\mu}$, each of the $ \alpha_{\mu
\lambda}^{(x_1,\ldots,x_r)}$'s is non-zero (these are indexed by the
$(k_1m_1, \ldots , k_rm_r)d$-paths from $\mu$ to $\lambda$).

\medskip
If $ \alpha_{\mu \lambda}^{(x_1,\ldots,x_r)} \neq 0$, then, since
there is a $(k_1m_1, \ldots , k_rm_r)d$-path from $\mu$ to
$\lambda$, and since the removal of a hook of length $md$ can be
obtained by the removal of a sequence of $m$ hooks of length $d$, we
see that $\mu$ and $\lambda$ have the same $d$-core.

\smallskip
We call the $ \alpha_{\mu \lambda}^{(x_1,\ldots,x_r)}$'s the
{\emph{MN-coefficients}}, and we will from now on write $
\alpha_{\mu \lambda}^{x}$ for $ \alpha_{\mu
\lambda}^{(x_1,\ldots,x_r)}$.

\subsection{One direction of the Nakayama Conjecture}

We are now able to prove that the unipotent $d$-blocks of $G$
satisfy one direction of a generalized Nakayama Conjecture. The
proof we give is an adaptation to our case of the proof given by
K\"{u}lshammer, Olsson and Robinson in the case of symmetric groups.

Take $x \in {\cal X}_d$ of $d$-type ${\bf km}=(k_1m_1, \ldots
,k_rm_r)$. We let $l=n-{\bf km}d=l-(\sum_{i=1}^{r} k_im_i)d$.
Writing $G_0=GL({\bf km}d,q)=\prod_{i=1}^{r} GL(k_im_id,q)$ and
$G_1=GL(l,q)$, we have $x =(x_0,x_1)=(x_0,1) \in G_0 \times G_1$.
Then $C_G(x)=H_0 \times H_1$, where $H_1 \cong GL(l,q)$ and
$H_0=C_{G_0}(x_0)$ (note that, if $x$ is semisimple, then $H_0 \cong
GL({\bf k},q^{{\bf m}d}) \cong \prod_{i=1}^r GL(k_i,q^{m_id})$).

Now take $y \in {\cal Y}_d(x)$. Then, as an element of $C_G(x)=H_0
\times H_1$, we have $y=(y_0,y_1)=(1,y_1)$.

Writing ${\cal Y}_d^t(u)$ for ${\cal Y}_d(u)$ when $u \in GL(t,q)$,
we have that $y$, element of $C_G(x)$, belongs to ${\cal Y}_d^n(x)$
if and only if $y=(1,y_1) \in H_0 \times H_1$, where $y_1$ belongs
to ${\cal Y}_d^l(1)$. Hence ${\cal Y}_d^n(x)$ is in natural one to
one correspondence with ${\cal Y}_d^l(1)$.
\medskip
Now we consider $\mu, \mu' \vdash n$, and $x \in {\cal X}_d$ of
$d$-type ${\bf km}=(k_1m_1, \ldots,k_rm_r)$. We have
$$\begin{array}{l} <\chi^{\mu}, \, \chi^{\mu'}>_{x{\cal Y}_d(x)} \displaystyle
 =\frac{1}{|G|} \sum_{y \in {\cal Y}_d(x)} \chi^{\mu}(xy) \overline{\chi^{\mu'}(xy)} \\
 =\displaystyle \frac{1}{|G|} \sum_{ y=(y_0,y_1) \in {\cal Y}_d(x)}  \left[ \left( \sum_{\lambda \in
 {\cal L}_{{\bf km}}^{\mu}} \alpha_{\mu \lambda}^x \chi^{\lambda}(y_1) \right)
 \left( \sum_{\lambda' \in {\cal L}_{{\bf km}}^{\mu'}} \alpha_{\mu' \lambda'}^x
 \overline{\chi^{\lambda'}(y_1)} \right) \right] \\  =\displaystyle \frac{1}{|G|} \sum_{y_1 \in
 {\cal Y}_d^l(1)} \; \;  \sum_{ \lambda \in {\cal L}_{{\bf km}}^{\mu} , \, \lambda' \in
 {\cal L}_{{\bf km}}^{\mu'}  } \alpha_{\mu \lambda}^x \alpha_{\mu' \lambda'}^x \chi^{\lambda}(y_1)
  \overline{\chi^{\lambda'}(y_1)} \end{array}$$
(by the above remark on ${\cal Y}_d(x)$).

\smallskip
We write $A_{\mu \mu'}^x=<\chi^{\mu}, \, \chi^{\mu'}>_{x {\cal
Y}(x)}$. Then
$$\begin{array}{cl} A_{\mu \mu'}^x &  =\displaystyle \frac{1}{|G|} \; \;  \sum_{  \lambda
\in {\cal L}_{{\bf km}}^{\mu}  , \, \lambda' \in {\cal L}_{{\bf km}}^{\mu'}  }
\alpha_{\mu \lambda}^x \alpha_{\mu' \lambda'}^x \sum_{y_1 \in {\cal Y}_d^l(1)} \chi^{\lambda}(y_1)
 \overline{\chi^{\lambda'}(y_1)} \\ &  =\displaystyle \frac{|H_1|}{|G|} \; \; \sum_{
 \lambda \in {\cal L}_{{\bf km}}^{\mu} , \, \lambda' \in {\cal L}_{{\bf km}}^{\mu'}  }
 \alpha_{\mu \lambda}^x \alpha_{\mu' \lambda'}^x <\chi^{\lambda} \chi^{\lambda'}>_{{\cal Y}_d^l(1)}
 \end{array}$$i.e.
$$ A_{\mu \mu'}^x =\displaystyle \frac{|H_1|}{|G|} \; \;
 \sum_{ \lambda \in {\cal L}_{{\bf km}}^{\mu} , \, \lambda' \in {\cal L}_{{\bf km}}^{\mu'} }
 \alpha_{\mu \lambda}^x \alpha_{\mu' \lambda'}^x A_{\lambda \lambda'}^1.$$
We use induction on $n$ to prove that, if $A_{\mu \mu'}^x \neq 0$,
then $\mu$ and $\mu'$ have the same $d$-core. We may assume that
$\mu \neq \mu'$.

If $n<d$, then each partition is its own $d$-core. Furthermore, in
this case, ${\cal X}_d=\{1\}$ and ${\cal Y}_d(1)=G$. Thus, for all
$x \in {\cal X}_d$, we have $A_{\mu \mu'}^x=A_{\mu \mu'}^1$, and
$A_{\mu \mu'}^1=<\chi^{\mu},\chi^{\mu'}>_G=0$ (since $\mu \neq
\mu'$). Hence the result is true in this case. Thus, we suppose $n
\geq d$. First suppose $x \neq 1$ and $A_{\mu \mu'}^x \neq 0$. Then
$$ A_{\mu \mu'}^x =\displaystyle \frac{|H_1|}{|G|} \; \; \sum_{ \lambda \in {\cal L}_{{\bf km}}^{\mu} , \, \lambda' \in {\cal L}_{{\bf km}}^{\mu'} } \alpha_{\mu \lambda}^x \alpha_{\mu' \lambda'}^x A_{\lambda \lambda'}^1 \neq 0.$$
Thus there exist $\lambda \in {\cal L}_{{\bf km}}^{\mu}$ and
$\lambda' \in {\cal L}_{{\bf km}}^{\mu'}$ such that  $\alpha_{\mu
\lambda}^x \alpha_{\mu' \lambda'}^x A_{\lambda \lambda'}^1 \neq 0$.
Then  $\alpha_{\mu \lambda}^x \neq 0$ implies that $\mu$ and
$\lambda$ have the same $d$-core, and $ \alpha_{\mu' \lambda'}^x
\neq 0$ implies that $\mu'$ and $\lambda'$ have the same $d$-core.
And, by the induction hypothesis (applied to $n-{\bf km}d <n$),
$A_{\lambda \lambda '}^1 \neq 0$ implies that $\lambda$ and
$\lambda'$ have the same $d$-core.

Now, if $x=1$, we see, by the existence and uniqueness of the
decomposition we introduced, that
$$\displaystyle 0=<\chi^{\mu}, \, \chi^{\mu'}>_G   = \sum_{ x \in {\cal X}_d }
<\chi^{\mu}, \, \chi^{\mu'}>_{x{\cal Y}_d(x)}  = \sum_{ x \in {\cal
X}_d  } A_{\mu \mu'}^x . $$ Hence, if $A_{\mu \mu'}^1 \neq 0$, then
there exists an $x' \in {\cal X}_d \setminus \{ 1 \}$ such that
$A_{\mu \mu'}^{x'} \neq 0$. This in turn implies, by the previous
case, that $\mu$ and $\mu'$ have the same $d$-core. Skipping back
from class functions to irreducible characters, we see that we have
proved the following

\begin{thm} Let $d>0$ be any integer. If two unipotent (irreducible) characters
$\chi_{\mu}$ and $\chi_{\mu'}$ of $G=GL(V)$ are directly linked
across some $x {\cal Y}_d(x)$, where $x \in {\cal X}_d$, then $\mu$
and $\mu'$ have the same $d$-core (and this is true in particular
for $x=1$).
\end{thm}

Extending by transitivity the relation of direct ${\cal
Y}_d(1)$-linking, we obtain
\begin{thm} If two unipotent characters $\chi_{\mu}$ and
$\chi_{\mu'}$ of $G$ are in the same unipotent $d$-block of $G$,
then $\mu$ and $\mu'$ have the same $d$-core.
\end{thm}

Each unipotent $d$-block of $G$ is therefore associated to a
$d$-core. For each given $d$-core $\gamma$, we can consider the
union of the (possibly many) unipotent $d$-blocks associated to
$\gamma$. The (a priori) bigger blocks obtained in this way are
parametrized by the set of $d$-cores of partitions of $n$, and they
satisfy an analogue of the Nakayama Conjecture. In accordance with
the terminology used in \cite{KOR}, we call them \emph{combinatorial
unipotent $d$-blocks}.

We would like to prove that, in fact, the notions of unipotent
$d$-block and combinatorial unipotent $d$-block coincide.
Unfortunately, it seems hard to prove (of infirm) in general. In the
following sections, we give examples where the computations are
easier, so that the result can be established.

\subsection{The case $d=1$}

In this section, we consider the extreme case $d=1$. Recall that the
1-regular elements of $G=GL(V)$ are just the unipotent elements.
Note that all partitions of $n$ have the same 1-core (the only
1-core there is, that is the empty one). In order to prove an
analogue of the Nakayama Conjecture for unipotent 1-blocks, we
therefore need to prove that any two unipotent characters belong to
the same 1-block. We actually prove more than that. Namely, any
irreducible character of $G$ is directly linked across unipotent
elements to the trivial character.

Our tool to deal with this case is the Alvis-Curtis duality. D.
Alvis has proved (cf \cite{Alvis}) that there is an isometry of
order 2, $D_G \colon ch(G) \longrightarrow ch(G)$, where $ch(G)$
denotes the ring of virtual characters of $G$ (in fact, this is true
in a much more general context, that is when $G$ is any finite group
with a split $(B,N)$-pair of characteristic $p$). In particular, the
image under the duality map $D_G$ of any irreducible (complex)
character of $G$ is, up to a sign, an irreducible character. The
case of unipotent characters has been made more precise by C. W.
Curtis (cf \cite{Curtis}). Namely, if $\chi_{\lambda} \in
$Unip$(G)$, then $D_G(\chi_{\lambda})=\chi_{\varepsilon \lambda}$,
where $\varepsilon$ is the sign character of the symmetric group
$S_n$ (that is, $\chi_{\varepsilon \lambda}$ is the unipotent
character $\chi_{\mu}$ where $\varphi_{\mu}$ is the irreducible
character of $S_n$ such that $\varphi_{\mu}=\varepsilon
\varphi_{\lambda}$). This can be written as
$D_G(\chi_{\lambda})=\chi_{\lambda^*}$, where $\lambda^*$ is the
partition conjugate to $\lambda$. In particular, we have
$D_G(1_G)=\chi_{\varepsilon}=St_G$, the Steinberg character of $G$.

\medskip
We define the virtual character $\chi_u$ of $G$ via
$$\chi_u(g)= \left\{ \begin{array}{l} 1 \; \mbox{if} \; g \;
\mbox{is unipotent} \\ 0 \; \mbox{otherwise} \end{array} \right. .$$
We then have the following:

\begin{lem}(\cite{Digne-Michel}, Lemme 11.1)
We have $\chi_u=D_G( |G|_{p'}^{-1} reg_G)$, where $reg_G$ is the
character of the regular representation of $G$.
\end{lem}

We let $G_u$ be the set of unipotent elements of $G$. Then, for any
$\chi \in$Irr$(G)$, we have
$$<\chi,1_G>_{G_u}:=\frac{1}{|G|} \sum_{g \in G_u} \chi(g)=\frac{1}{|G|} \sum_{g \in G}
\chi(g) \chi_u(g)=<\chi,\chi_u>_G.$$ Thus, since $D_G$ is an
isometry, we have
$$<\chi,1_G>_{G_u}=<D_G(\chi),D_G(\chi_u)>_G=\frac{1}{|G|_{p'}}
<D_G(\chi),reg_G>_G=\frac{1}{|G|_{p'}}D_G(\chi)(1).$$ Hence
$<\chi,1_G>_{G_u} \neq 0$, i.e. $\chi$ and $1_G$ are directly linked
across unipotent elements. We thus have

\begin{thm}
The only 1-block of $G$ is Irr$(G)$, and the only unipotent 1-block
of $G$ is Unip$(G)$. In particular, the unipotent 1-blocks of $G$
satisfy an analogue of the Nakayama Conjecture.
\end{thm}

Note that, if $\chi=\chi_{\lambda} \in$Unip$(G)$, then we have
$<\chi_{\lambda},1_G>_{G_u}=\frac{1}{|G|_{p'}}\chi_{\lambda^*}(1)$.

\subsection{$d$-weight and simple partitions}

From now on, we fix the integers $n$ and $d$, $d \neq 1$, and some
$d$-core $\gamma$. Our aim is to prove that, if two partitions
$\lambda$ and $\mu$ of $n$ have the same $d$-core $\gamma$, then the
unipotent characters $\chi_{\lambda}$ and $\chi_{\mu}$ belong to the
same unipotent $d$-block of $G=GL(V)$. However, if we simply write
down the inner product of $\chi_{\lambda}$ and $\chi_{\mu}$ on
$d$-regular (or $d$-singular) elements of $G$, we obtain a huge
expression which is very hard to manipulate, let alone trying to
prove that it is non-zero.

Instead, what we do in sections 4.5 and 4.6 is build, under some
extra assumption on $\lambda$ and $\mu$, a chain of unipotent
characters $\chi_0=\chi_{\lambda}, \, \chi_1, \, \ldots, \,
\chi_r=\chi_{\mu}$, such that, for each $0 \leq i \leq r-1$,
$\chi_i$ and $\chi_{i+1}$ are directly linked across $d$-singular
elements. This proves that $\chi_{\lambda}$ and $\chi_{\mu}$ are
linked across $d$-singular elements, and therefore belong to the
same unipotent $d$-block.

Section 4.5 is devoted to the statement and proof of Theorem 4.6,
which corresponds to the "elementary links" in our chain of
characters. In section 4.6, we then show that, if $\lambda$ and
$\mu$ have "weight small compared to $d$", then it is possible to
build such a chain between $\chi_{\lambda}$ and $\chi_{\mu}$ (cf
Theorem 4.10).

\smallskip
We now introduce the observations and notions we will use to
construct our chain.

We can associate to each combinatorial unipotent $d$-block a
$d${\emph{-weight}}, the $d$-weight of any partition of $n$ labeling
some unipotent character in the block.

The first observation comes from the Murnaghan-Nakayama Rule: if $\lambda$ has $d$-weight
 $w$, then $\chi_{\lambda}$ vanishes on any element of $G$ whose $d$-part has $d$-type
  $(k_1 m_1, \, \ldots , \, k_r m_r)$ with $\sum_i k_i m_i >w$.

This leads to the first definition and simplification. For any
element $g \in G$ with  $d$-part of $d$-type $(k_1 m_1, \, \ldots ,
\, k_r m_r)$, we define the $d${\emph{-weight}} of $g$ to be $\sum_i
k_i m_i$. Thus, for partitions of $d$-weight $w$, it is enough to
consider the inner product of $\chi_{\lambda}$ and $\chi_{\mu}$ on
$d$-singular elements of $d$-weight {\bf{at most}} $w$.

\smallskip

The Murnaghan-Nakayama Rule also implies that, if, for all $m>1$, ${\lambda}$ contains no
$md$-hook, then $\chi_{\lambda}$ vanishes on any element $g$ with $d$-part of $d$-type
$(k_1 m_1, \, \ldots , \, k_r m_r)$ when there exists $i$ such that $m_i>1$ (i.e. whenever
 there exists an irreducible factor of degree $md$ with $m>1$ in $Min(g)$). If $\lambda$ is
  such a partition, we will say that $\lambda$ is {\emph{simple}} (or $d${\emph{-simple}}).

\medskip
Given any partition $\lambda$ of $n$, a useful and convenient way to
store the $d$-information about $\lambda$ is with an abacus. For a
complete definition, we refer to \cite{James-Kerber}, Section 2.7
(note however that the abacus we describe here is the horizontal
mirror image of that used by James and Kerber).

We present the construction of the abacus with an example: we take
$d=3$ and the partition $\lambda=(6,5,5,2,1)$ of $n=19$. We put the
Young diagram of $\lambda$ in the upper-left corner of a quarter
plane, and consider its (infinite) rim, in bold below. We choose an
arbitrary origin, indicated by the symbol $\triangleright$, on the
vertical axis.

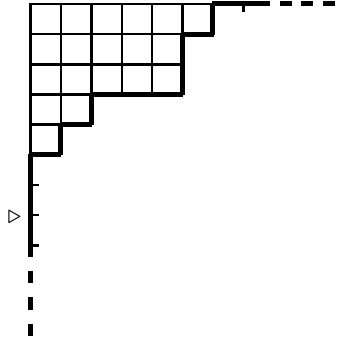
\setlength{\unitlength}{1mm}

\begin{center}

\begin{picture}(44,44)

\linethickness{0.5mm}

\put(0,24) {\line(1,0){4}}

\put(4,28) {\line(1,0){4}}

\put(8,32) {\line(1,0){12}}

\put(20,40){\line(1,0){4}}

\put(24,44){\line(1,0){6}}

\dashline[+60]{2}(30,44)(40,44)

\put(0,12){\line(0,1){12}}

\put(-3,15){$\triangleright$}

 \dashline[+70]{2}(0,12)(0,0)

\put(4,24){\line(0,1){4}}

\put(8,28) {\line(0,1){4}}

\put(20,32) {\line(0,1){8}}

\put(24,40){\line(0,1){4}}

 \linethickness{0.2mm}

\put(0,28) {\line(1,0){4}}

\put(0,32) {\line(1,0){8}}

\put(0,36) {\line(1,0){20}}

\put(0,40) {\line(1,0){20}}

\put(0,44) {\line(1,0){24}}

\put(0,24) {\line(0,1){20}}

\put(4,28) {\line(0,1){16}}

\put(8,32) {\line(0,1){12}}

\put(12,32) {\line(0,1){12}}

\put(16,32) {\line(0,1){12}}

\put(20,40) {\line(0,1){4}}

\put(0,16){\line(1,0){1}}

\put(0,12){\line(1,0){1}}

\put(0,20){\line(1,0){1}}

\put(28,43){\line(0,1){1}}

\end{picture}

\end{center}

We see the rim as an infinite sequence of vertical and horizontal
dashes of length 1 (the length of a box in the diagram). Writing 1
for a vertical dash and 0 for a horizontal one, we encode the rim as
the following sequence, where we indicate the origin as before:
$$ \cdots \; \; 1 \; \; 1  \; \; \triangleright \; \; 1 \; \; 1 \; \; 0 \; \; 1 \; \; 0 \; \; 1 \; \;
0 \; \; 0 \; \; 0 \; \; 1 \; \; 1 \; \; 0 \; \; 1  \; \; 0 \; \; 0
\; \; 0 \; \; \cdots  \; \; \; \; \; \; \; (\dagger)$$ We now
construct the abacus using this sequence. We put beads on $d=3$
runners, going from left to right and bottom to top, putting a bead
for each 1 and an empty spot for each 0. Here again, we indicate the
origin. We get

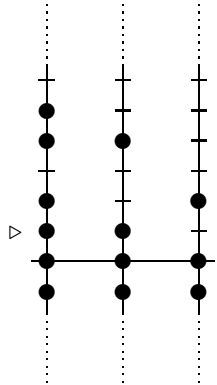
\setlength{\unitlength}{1mm}

\begin{center}

\begin{picture}(24,50)

\put(2,9) {\line(0,1){33}}

\put(12,9) {\line(0,1){33}}

\put(22,9) {\line(0,1){33}}

\put(0,16) {\line(1,0){24}}

\put(2,12){\circle*{2}}

\put(12,12){\circle*{2}}

\put(22,12){\circle*{2}}

\put(2,16){\circle*{2}}

\put(12,16){\circle*{2}}

\put(22,16){\circle*{2}}

\put(2,20){\circle*{2}}

\put(12,20){\circle*{2}}

\put(2,24){\circle*{2}}

\put(22,24){\circle*{2}}

\put(2,32){\circle*{2}}

\put(12,32){\circle*{2}}

\put(2,36){\circle*{2}}

\put(1,28) {\line(1,0){2}}

\put(11,28) {\line(1,0){2}}

\put(21,28) {\line(1,0){2}}

\put(11,24) {\line(1,0){2}}

\put(21,20) {\line(1,0){2}}

\put(11,36) {\line(1,0){2}}

\put(21,36) {\line(1,0){2}}

\put(21,32) {\line(1,0){2}}

\put(1,40) {\line(1,0){2}}

\put(11,40) {\line(1,0){2}}

\put(21,40) {\line(1,0){2}}

\dottedline{1}(2,42)(2,50)

\dottedline{1}(12,42)(12,50)

\dottedline{1}(22,42)(22,50)

\dottedline{1}(12,0)(12,9)

\dottedline{1}(2,0)(2,9)

\dottedline{1}(22,0)(22,9)

\put(-3,19){$\triangleright$}

\end{picture}

\end{center}

We see that, up to the choice of the runner where we put the origin,
there is bijection that associates to each partition a unique
abacus.

The $d$-information is visible in the abacus in a natural way. If $k$ is a positive integer, then any $kd$-rim-hook in $\lambda$ corresponds to a sequence of $kd+1$ digits in $(\dagger)$ which starts with a 0 and ends with a 1, and to a bead in the abacus which lies, on the same runner, $k$ places above an empty spot. Moreover, these correspondences are bijective.

The removal of a $kd$-rim-hook in the Young diagram is achieved by exchanging the 0 and 1 at the extremity of the corresponding subsequence of $(\dagger)$, and by moving the corresponding bead to the empty spot.

By removing all the $d$-hooks from its Young diagram, we get the $d$-core of $\lambda$, and the corresponding sequence and abacus:

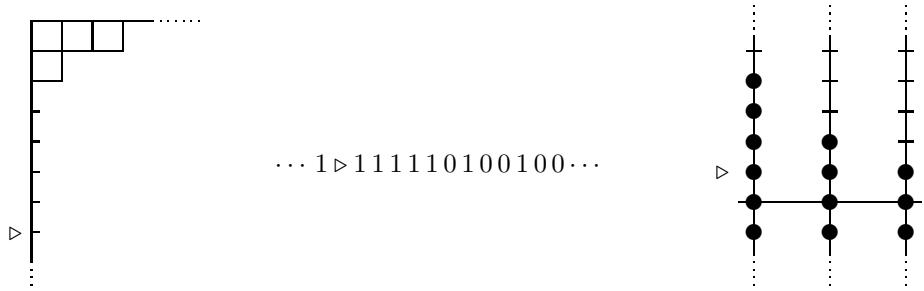
\begin{center}

\begin{picture}(115,42)

\put(0,8) {\line(0,1){32}}

\put(-3,11) {$\triangleright$}

\dottedline{1}(0,8)(0,4)

\put(0,40) {\line(1,0){16}}

\dottedline{1}(16,40)(22,40)

\put(0,36) {\line(1,0){12}}

\put(0,32) {\line(1,0){4}}

\put(4,32){\line(0,1){8}}

\put(8,36){\line(0,1){4}}

\put(12,36){\line(0,1){4}}

\put(0,12) {\line(1,0){1}}

\put(0,16) {\line(1,0){1}}

\put(0,20) {\line(1,0){1}}

\put(0,24) {\line(1,0){1}}

\put(0,28) {\line(1,0){1}}

\put(32,20){$\cdots \, 1 \triangleright 1 \, 1 \, 1 \, 1 \, 1 \, 0
\, 1 \, 0 \, 0 \, 1 \, 0 \, 0 \cdots$}

\put(95,9) {\line(0,1){29}}

\put(105,9) {\line(0,1){29}}

\put(115,9) {\line(0,1){29}}

\put(93,16) {\line(1,0){24}}

\put(95,12){\circle*{2}}

\put(105,12){\circle*{2}}

\put(115,12){\circle*{2}}

\put(95,16){\circle*{2}}

\put(105,16){\circle*{2}}

\put(115,16){\circle*{2}}

\put(95,20){\circle*{2}}

\put(105,20){\circle*{2}}

\put(115,20){\circle*{2}}

\put(95,24){\circle*{2}}

\put(105,24){\circle*{2}}

\put(95,32){\circle*{2}}

\put(95,28){\circle*{2}}

\put(94,36) {\line(1,0){2}}

\put(104,28) {\line(1,0){2}}

\put(114,28) {\line(1,0){2}}

\put(114,24) {\line(1,0){2}}

\put(114,36) {\line(1,0){2}}

\put(104,36) {\line(1,0){2}}

\put(114,32) {\line(1,0){2}}

\put(104,32) {\line(1,0){2}}

\dottedline{1}(95,38)(95,42)

\dottedline{1}(105,38)(105,42)

\dottedline{1}(115,42)(115,38)

\dottedline{1}(95,4)(95,9)

\dottedline{1}(105,4)(105,9)

\dottedline{1}(115,4)(115,9)

\put(90,19){$\triangleright$}

\end{picture}

\end{center}

To each runner of the abacus, we associate a partition as follows. When we move all the beads as far down as possible, we get one part for each bead that we move, of length the number of places the bead goes down. The resulting $d$-tuple $(\lambda^{(1)}, \, \ldots , \, \lambda^{(d)})$ of
partitions is the $d$-quotient of $\lambda$. In the above example, we get $(\lambda^{(1)}, \, \lambda^{(2)} , \, \lambda^{(3)})=((1,1),(2),(1))$. The lengths of the $\lambda^{(i)}$'s add up
to the $d$-weight $w$ of $\lambda$. For any positive integer $k$, we call $k$-hook in the quotient of $\lambda$ any $k$-hook in one of the $\lambda^{(i)}$'s. Then there is a bijection between the set of $k$-hooks in the $d$-quotient of $\lambda$ and the set of $kd$-hooks in $\lambda$.

One sees easily that the partition $\lambda$ is uniquely determined by its $d$-core and
$d$-quotient.

\medskip

From the definition we gave above and the description of the abacus,
we see that, if $\lambda$ is $d$-simple of weight $w$, then the
moves corresponding to the removals of all the $d$-hooks from
$\lambda$ must be stored on exactly $w$ distinct runners of the
abacus (so that, in particular, one must have $w \leq d$), and they
are all elementary, i.e. involve each a single bead which can be
taken exactly once one place down.

\smallskip

We say that two partitions are {\emph{disjoint}} if the possible
moves on their abacus are stored on different runners. In terms of
quotients, we see that two partitions $\lambda$ and $\mu$, with
$d$-quotients $(\lambda^{(1)}, \, \ldots , \, \lambda^{(d)})$ and
$(\mu^{(1)}, \, \ldots , \, \mu^{(d)})$ respectively, are disjoint
if and only if, for each $1 \leq i \leq d$, $\lambda^{(i)} \ne
\emptyset $ implies  $\mu^{(i)}=\emptyset$.

For any partition $\lambda$ and integer $i \geq 0$, we write ${\cal
L}_{\lambda}^i$ (respectively ${\cal L}_{\lambda}^{(i)}$) the set of
partitions we can get by removing $i$ $d$-hooks (respectively one
$id$-hook) from $\lambda$. In particular, ${\cal L}_{\lambda}^0=\{
\lambda \}$, and, if $\lambda$ has $d$-weight $w$ and $d$-core
$\gamma$, then ${\cal L}_{\lambda}^w=\{ \gamma \}$ and, for all $i >
w$, ${\cal L}_{\lambda}^i=\emptyset$.

\subsection{One result of direct linking}

Suppose that $\lambda$ and $\mu$ are partitions of $n$ of $d$-weight $w \neq 0$, and with the same $d$-core $\gamma$. Suppose furthermore that $\mu$ is simple and disjoint from $\lambda$. Then, by definition(s), we have:

\begin{itemize}
\item{}${\cal L}_{\mu}^{(i)}= \emptyset$ for $1<i \leq w$.
\item{} ${\cal L}_{\lambda}^i \cap {\cal L}_{\mu}^i= \emptyset$ for $0\leq i<w$.\item{} ${\cal L}_{\lambda}^w={\cal L}_{\mu}^w= \{ \gamma \}$.
\end{itemize}

These properties induce a large number of simplifications in the expression for the inner product of $\chi_{\lambda}$ and $\chi_{\mu}$ on $d$-singular elements.
We will prove the following:

\begin{thm}
Suppose $\lambda$ and $\mu$ are partitions of $n$ of $d$-weight $w
\neq 0$, and with the same $d$-core $\gamma$. Suppose furthermore
that $\mu$ is simple and disjoint from $\lambda$. Then the unipotent
characters $\chi_{\lambda}$ and $\chi_{\mu}$ are directly linked
across the set of $d$-regular elements of $G$.
\end{thm}

We write $F=|{\cal F}_d|$, where ${\cal F}_d$ is the set of
irreducible monic polynomials of degree $d$ over ${\bf F}_q$, and we
suppose that $F \geq n/d$ (so that $F \geq w$ for any possible
$d$-weight $w$).

We use induction on $w$ to prove that, if $\lambda$ and $\mu$ are as
above, then the inner product of $\chi^{\lambda}$ and $\chi^{\mu}$
on the set of $d$-regular elements of $G$ is$$<\chi^{\lambda}, \,
\chi^{\mu}>_{d-reg}= \displaystyle \frac{ (-1)^w F^w}{w! (q^d-1)^w}
\; \; \sum_{P_{\lambda} \in {\cal P}^{\lambda}_{\gamma}, \, P_{\mu}
\in {\cal P}^{\mu}_{\gamma}} (-1)^{L_{P_{\lambda}}}
(-1)^{L_{P_{\mu}}},$$ where ${\cal P}^{\lambda}_{\gamma}$ (resp.
${\cal P}^{\mu}_{\gamma}$) denotes the set of paths from $\lambda$
(resp. $\mu$) to $\gamma$, in the lattice of partitions, obtained by
successive removals of $d$-hooks, and, for $P_{\lambda} \in {\cal
P}^{\lambda}_{\gamma}$ (resp. $P_{\mu} \in {\cal
P}^{\mu}_{\gamma}$), $L_{P_{\lambda}}$ (resp. $L_{P_{\mu}}$) is the
sum of the leg lengths of the hooks removed along the corresponding
path.

\medskip

We suppose the result established for partitions of $d$-weight up to
$w-1 \neq 0$, and we take $\lambda$ and $\mu$ of weight $w$ as
above. We order the elements of ${\cal F}_d$, and write ${\cal F}_d=
\{ f_1, \, \ldots , \, f_F \}$. We write ${\cal C}$ for the set of
$d$-singular elements of weight at most $w$ of $G$. The conjugacy
classes in ${\cal C}$ are parametrized by the $F$-tuples of
(possibly empty) partitions $(\nu_1, \, \ldots , \, \nu_F)$, where
$\nu_i \vdash k_i$ for $1 \leq i \leq F$ and $1 \leq k_1 + \cdots +
k_F \leq w$. We choose a basis for $V$ and identify $G=GL(V)$ with
$GL(n,q)$. We will write $G_n$ for $GL(n,q)$, and similarly for
smaller dimensions. Writing $\sim$ for conjugation, we get
$${\cal C}/ \sim  \; = \displaystyle \bigcup_{k=1}^w \left\{ \left( \begin{array}{cc} \rho & \\  & \sigma \end{array} \right) , \, \rho \in {\cal C}_k/ \sim, \, \sigma \in G_{n-kd} , \, \sigma \;  d-reg \right\},$$
where
$${\cal C}_k/ \sim \; = \displaystyle \bigcup_{k_1 + \cdots +k_F=k} \; \; \bigcup_{\nu_1 \vdash k_1, \,
\ldots, \, \nu_F \vdash k_F} \left\{ \left( \begin{array}{ccc}
U_{\nu_1}(f_1) & & \\  & \ddots & \\ & & U_{\nu_F}(f_F) \end{array}
\right) \right\}.$$

We have
$$ <\chi^{\lambda}, \, \chi^{\mu}>_{d-sing} = <\chi^{\lambda}, \, \chi^{\mu}>_{\cal C} = \displaystyle \sum_{k=1}^w <\chi^{\lambda}, \, \chi^{\mu}>_{{\cal C}_k} $$
and, for $1 \leq k \leq w$, $<\chi^{\lambda}, \, \chi^{\mu}>_{{\cal C}_k}$ is equal to
$$ \displaystyle  \sum_{k_1 + \cdots +k_F=k} \; \; \sum_{\nu_1 \vdash k_1, \, \ldots, \, \nu_F \vdash k_F} \; \; \sum_{\scriptstyle \sigma \in G_{n-kd} \atop \scriptstyle \sigma \; d-reg} \frac{\chi^{\lambda}(\rho^{\nu_1} \ldots \rho^{\nu_F} \sigma)\overline{\chi^{\mu}(\rho^{\nu_1} \ldots \rho^{\nu_F} \sigma)}}{|C_n(\rho^{\nu_1} \ldots \rho^{\nu_F} \sigma)|},$$
where we write $C_n(g)$ for $C_{G_n}(g)=C_{GL(n,q)}(g)$ and
$\rho^{\nu}$ for $U_{\nu}(f)$, $f \in {\cal F}$ (abuse justified
since the values of unipotent characters on $U_{\nu}(f)$ and the
size of its centralizer don't depend on $f$, but only on $\nu$).

Repeated use of the Murnaghan-Nakayama Rule gives us
$$\chi^{\lambda}(\rho^{\nu_1} \ldots \rho^{\nu_F} \sigma)= \; \; \displaystyle \sum_{\lambda_1 \in {\cal L}_{\lambda}^{k_1}, \, \lambda_2 \in {\cal L}_{\lambda_1}^{k_2}, \, \ldots , \, \lambda_F \in {\cal L}_{\lambda_{F-1}}^{k_F} } a_{\lambda \lambda_1}^{\nu_1}  a_{\lambda_1 \lambda_2}^{\nu_2} \ldots a_{\lambda_{F-1} \lambda_F}^{\nu_F} \chi^{\lambda_F}(\sigma),$$
and a similar expression for $\chi^{\mu}(\rho^{\nu_1} \ldots \rho^{\nu_F} \sigma)$.

We obtain
$$ <\chi^{\lambda}, \, \chi^{\mu}>_{d-sing}=
 \displaystyle \sum_{k=1}^w \; \; \sum_{k_1 + \cdots +k_F=k}  {\cal A}_{k_1,\ldots,k_F}^{\lambda,\mu},$$
where $${\cal A}_{k_1,\ldots,k_F}^{\lambda,\mu}=\displaystyle \sum_{ \lambda_1 \in {\cal L}_{\lambda}^{k_1}, \, \ldots , \, \lambda_F \in {\cal L}_{\lambda_{F-1}}^{k_F}} \; \; \sum_{\mu_1 \in {\cal L}_{\mu}^{k_1}, \,  \ldots , \, \mu_F \in {\cal L}_{\mu_{F-1}}^{k_F}} {\cal B}_{\lambda_1,\mu_1}^{\lambda,\mu} \ldots {\cal B}_{\lambda_F,\mu_F}^{\lambda_{F-1},\mu_{F-1}} {\cal D}_{\lambda_F,\mu_F},$$
with
$${\cal B}_{\lambda_1,\mu_1}^{\lambda,\mu}=\displaystyle \sum_{\nu_1 \vdash k_1} \frac{a_{\lambda \lambda_1}^{\nu_1}a_{\mu \mu_1}^{\nu_1}}{|C_{k_1d}(\rho^{\nu_1})|} , \;  \ldots , \;   {\cal B}_{\lambda_F,\mu_F}^{\lambda_{F-1},\mu_{F-1}}= \sum_{\nu_F \vdash k_F} \frac{a_{\lambda_{F-1} \lambda_F}^{\nu_F}a_{\mu_{F-1} \mu_F}^{\nu_F}}{|C_{k_Fd}(\rho^{\nu_F})|}$$
and
$$ {\cal D}_{\lambda_F,\mu_F}= \displaystyle \sum_{\scriptstyle \sigma \in G_{n-kd}/ \sim \atop \scriptstyle
\sigma \; d-reg} \frac{\chi^{\lambda_F}(\sigma)
\overline{\chi^{\mu_F}(\sigma)}}{|C_{n-kd}(\sigma)|}.$$ Now, for all
$1 \leq k \leq w$, we have
$$\displaystyle \sum_{\scriptstyle \sigma \in G_{n-kd}/ \sim \atop \scriptstyle \sigma \; d-reg}
 \frac{\chi^{\lambda_F}(\sigma) \overline{\chi^{\mu_F}(\sigma)}}{|C_{n-kd}(\sigma)|}=<\chi^{\lambda_F}, \, \chi^{\mu_F}>_{G_{n-kd}, \, d-reg},$$
which is known, by the induction hypothesis, since $\lambda_F$ and $\mu_F$ have weight $w-k<w$ (and the same $d$-core $\gamma$).

In order to deal with the other factors, we start with the following lemma, which doesn't depend on the special choice of $\mu$:

\begin{lem}
For any $1 \leq i \leq F$, we have, writing $\lambda_0$ and $\mu_0$ for $\lambda$ and $\mu$ respectively,
$$\displaystyle \sum_{\nu_i \vdash k_i} \frac{a_{\lambda_{i-1} \lambda_i}^{\nu_i}a_{\mu_{i-1} \mu_i}^{\nu_i}}{|C_{k_id}(\rho^{\nu_i})|}=\varepsilon_{\lambda_{i-1}\lambda_i} \varepsilon_{\mu_{i-1} \mu_i} \sum_{\alpha \vdash k_i} \frac{\phi_{\lambda_{i-1} | \lambda_i}(\alpha) \phi_{\mu_{i-1} | \mu_i}(\alpha)}{|C_{S_{k_i}}(\alpha)| \; |T_{\overline{\alpha}}| },$$
where:
\begin{itemize}
\item{}
if $\alpha=(1^{r_1}, \, 2^{r_2}, \, \ldots , \, s^{r_s}) \vdash k_i$, then $\overline{\alpha}=(d^{r_1}, \, (2d)^{r_2}, \, \ldots , \, (sd)^{r_s}) \vdash k_id$,
\item
$T_{\overline{\alpha}}$ is a maximal torus of $G_{k_id}=GL(k_id,q)$
representing the $G_{k_id}$-conjugacy class corresponding to the
$S_{k_i}$-conjugacy class of $\alpha$;\\ we have
$|T_{\overline{\alpha}}|=\prod_i (q^{di}-1)^{r_i}$ (cf \cite{Fo-Sri}
(1.12));
\item
 $\phi_{\mu_{i-1} | \mu_i} \in {\bf Z} \mbox{Irr}(S_{k_i})$, and, if $h \in S_{\ell}$ and $g \in S_{k_id}$ is
  of cycle type $\overline{\alpha}$ ($\alpha \vdash k_i$), then
$$\varphi_{\mu_{i-1}}(gh)=\displaystyle \sum_{\eta \in {\cal L}_{\mu_{i-1}}^{k_i}}
 \varepsilon_{\mu_{i-1} \eta} \phi_{\mu_{i-1} | \eta}(\alpha) \varphi_{\eta}(h)$$
(where, for any positive integer $m$, we write Irr$(S_m)=\{
\varphi_{\zeta}, \, \zeta \vdash m \}$) and similarly for
$\phi_{\lambda_{i-1} | \lambda_i}$ (cf \cite{Farahat} and
\cite{Fo-Sri} (2.3)).
\end{itemize}
\end{lem}
\begin{proof}

For all $\nu_i \vdash k_i$, we have $C_{k_id}(\rho^{\nu_i})=C_{C_{k_id}(\rho^{\nu_i}_S)}(\rho^{\nu_i}_U)$, where $\rho^{\nu_i}_S$ and $\rho^{\nu_i}_U$ are the semisimple and unipotent parts of $\rho^{\nu_i}$ respectively, and $C_{k_id}(\rho^{\nu_i}_S) \cong GL(k_i, q^d)$. We write $H_i$ for $GL(k_i, q^d)$ (then $H_i$ has Weyl group $W_{H_i} \cong S_{k_i}$), so that $C_{k_id}(\rho^{\nu_i}) \cong C_{H_i}(\rho^{\nu_i}_U)$.

We have, writing $H_i^u/ \sim$ for a set of representatives for the unipotent conjugacy classes of $H_i$,
$$\displaystyle \sum_{\nu_i \vdash k_i} \frac{a_{\lambda_{i-1} \lambda_i}^{\nu_i}a_{\mu_{i-1} \mu_i}^{\nu_i}}{|C_{k_id}(\rho^{\nu_i})|}= \sum_{\rho^{\nu_i}_U \in H_i^u/ \sim } \frac{a_{\lambda_{i-1} \lambda_i}^{\nu_i}a_{\mu_{i-1} \mu_i}^{\nu_i}}{|C_{H_i}(\rho^{\nu_i}_U)|}.$$
Now, by definition (cf \cite{Fo-Sri} Theorem (2G)), we have
$$\displaystyle a_{\lambda_{i-1} \lambda_i}^{\nu_i}= \varepsilon_{\lambda_{i-1}\lambda_i} \sum_{\alpha \vdash k_i} \frac{1}{|W_{H_i}(T_{\overline{\alpha}})|} Q_{T_{\overline{\alpha}}}^{H_i}(\rho^{\nu_i}_U)\phi_{\lambda_{i-1} | \lambda_i}(\alpha)$$
(and similarly for $a_{\mu_{i-1} \mu_i}^{\nu_i}$), where $W_{H_i}(T_{\overline{\alpha}})$ is the Weyl group of $T_{\overline{\alpha}}$ in $H_i$, so that $|W_{H_i}(T_{\overline{\alpha}})|=|C_{S_{k_i}}(\alpha)|$ (cf \cite{Fo-Sri} (1.11)), and $Q_{T_{\overline{\alpha}}}^{H_i}$ is the Green function, integer-valued class function defined on $H_i^u$, associated to $T_{\overline{\alpha}}$ (cf \cite{Fo-Sri}).

We obtain that $\sum_{\nu_i \vdash k_i} \frac{a_{\lambda_{i-1} \lambda_i}^{\nu_i}a_{\mu_{i-1} \mu_i}^{\nu_i}}{|C_{k_id}(\rho^{\nu_i})|}$ is equal to
$$\sum_{\rho^{\nu_i}_U \in H_i^u/ \sim } \frac{\varepsilon_{\lambda_{i-1}\lambda_i} \varepsilon_{\mu_{i-1} \mu_i}}{|C_{H_i}(\rho^{\nu_i}_U)|}  \sum_{\alpha, \beta \vdash k_i} \frac{ Q_{T_{\overline{\alpha}}}^{H_i}(\rho^{\nu_i}_U)\phi_{\lambda_{i-1} | \lambda_i}(\alpha)  Q_{T_{\overline{\beta}}}^{H_i}(\rho^{\nu_i}_U) \phi_{\mu_{i-1} | \mu_i}(\beta)}{|W_{H_i}(T_{\overline{\alpha}})| \; |W_{H_i}(T_{\overline{\beta}})| } ,$$
which, since the Green functions are class functions, can be written as
$$\varepsilon_{\lambda_{i-1}\lambda_i} \varepsilon_{\mu_{i-1} \mu_i} \sum_{\alpha, \beta \vdash k_i} \frac{\phi_{\lambda_{i-1} | \lambda_i}(\alpha) \phi_{\mu_{i-1} | \mu_i}(\beta)}{|W_{H_i}(T_{\overline{\alpha}})| \; |W_{H_i}(T_{\overline{\beta}})| } \left( \frac{1}{|H_i|} \sum_{\rho^{\nu_i}_U \in H_i^u} Q_{T_{\overline{\alpha}}}^{H_i}(\rho^{\nu_i}_U)  Q_{T_{\overline{\beta}}}^{H_i}(\rho^{\nu_i}_U) \right).$$
Now, by the first orthogonality relation for Green functions in $H_i$ (cf \cite{Fo-Sri} (1.7)), we have
$$\frac{1}{|H_i|} \sum_{\rho^{\nu_i}_U \in H_i^u} Q_{T_{\overline{\alpha}}}^{H_i}(\rho^{\nu_i}_U)  Q_{T_{\overline{\beta}}}^{H_i}(\rho^{\nu_i}_U) = \left\{ \begin{array}{ll} \frac{|W_{H_i}(T_{\overline{\alpha}})|}{|T_{\overline{\alpha}}|} & \mbox{if} \; T_{\overline{\alpha}} \sim_{H_i} T_{\overline{\beta}} \\ 0 & \mbox{otherwise} \end{array} \right. .$$
But $ T_{\overline{\alpha}} \sim_{H_i} T_{\overline{\beta}}$ if and only if $\overline{\alpha}=\overline{\beta}$, if and only if $\alpha=\beta$. We thus get, since $|W_{H_i}(T_{\overline{\alpha}})|=|C_{S_{k_i}}(\alpha)|$,
$$\sum_{\nu_i \vdash k_i} \frac{a_{\lambda_{i-1} \lambda_i}^{\nu_i}a_{\mu_{i-1} \mu_i}^{\nu_i}}{|C_{k_id}(\rho^{\nu_i})|}=\varepsilon_{\lambda_{i-1}\lambda_i} \varepsilon_{\mu_{i-1} \mu_i} \sum_{\alpha \vdash k_i} \frac{\phi_{\lambda_{i-1} | \lambda_i}(\alpha) \phi_{\mu_{i-1} | \mu_i}(\alpha)}{|C_{S_{k_i}}(\alpha)| \; |T_{\overline{\alpha}}| },$$
which ends the proof.
\end{proof}

The fact that $\mu$ is simple then gives us further simplifications; we have the following
\begin{lem}
For any $1 \leq i \leq F$ and $\alpha \in S_{k_i}$, we have $\phi_{\mu_{i-1} | \mu_i}(\alpha)=0$ unless $\alpha=(1^{k_i})$ (the identity element of $S_{k_i}$).\end{lem}
\begin{proof}
Take any $1 \leq i \leq F$. We have $\mu_{i-1} \vdash m$ and $\mu_i
\vdash \ell$, with $m=\ell + k_i d$. Then, by definition (cf
\cite{Fo-Sri} (2.3) and \cite{Farahat}), $\phi_{\mu_{i-1} | \mu_i}
\in {\bf Z} \mbox{Irr}(S_{k_i})$, and, if $h \in S_{\ell}$ and $g
\in S_{k_id}$ is of cycle type $\overline{\alpha}$ ($\alpha \vdash
k_i$), then
$$\varphi_{\mu_{i-1}}(gh)=\displaystyle \sum_{\eta \in {\cal L}_{\mu_{i-1}}^{k_i}} \varepsilon_{\mu_{i-1} \eta} \phi_{\mu_{i-1} | \eta}(\alpha) \varphi_{\eta}(h).$$
Now $\mu_{i-1}$ can be obtained from $\mu$, which is $d$-simple, by removing a sequence of $d$-hooks, so that $\mu_{i-1}$ is itself simple, whence, for $j>1$, there is no $jd$-hook in $\mu_{i-1}$.

Thus, by the Murnaghan-Nakayama rule (in the symmetric group), if there is a cycle of length $j>1$ in $\alpha$ (so that there is a $jd$-cycle in $g$), then
$$\varphi_{\mu_{i-1}} (gh)=0, \; \mbox{for all} \; h \in S_{\ell}.$$
Hence, for all $h \in S_{\ell}$,
$$\displaystyle \sum_{\eta \in {\cal L}_{\mu_{i-1}}^{k_i}} \varepsilon_{\mu_{i-1} \eta} \phi_{\mu_{i-1} | \eta}(\alpha) \varphi_{\eta}(h)=0,$$
whence
$$ \varepsilon_{\mu_{i-1} \eta} \phi_{\mu_{i-1} | \eta}(\alpha)=0 \; \mbox{for all} \; \eta \in {\cal L}_{\mu_{i-1}}^{k_i},$$
in particular for $\eta = \mu_i$.
\end{proof}

By the previous two lemmas, we therefore get, for any $1 \leq i \leq F$,
$$ {\cal B}_{\lambda_i,\mu_i}^{\lambda_{i-1},\mu_{i-1}}= \displaystyle \frac{
 \varepsilon_{\lambda_{i-1} \lambda_i}  \varepsilon_{\mu_{i-1} \mu_i} \phi_{\lambda_{i-1} | \lambda_i}(1^{k_i}) \phi_{\mu_{i-1} | \mu_i}(1^{k_i}) }{|C_{S_{k_i}}(1^{k_i})| \, |T_{(\overline{1^{k_i}})}| },$$
and we have $C_{S_{k_i}}(1^{k_i})=S_{k_i}$ of order $k_i!$, and $ |T_{(\overline{1^{k_i}})}|= |T_{(d^{k_i})}|=(q^d-1)^{k_i}$.

Now, we have, for all $h \in S_{\ell}$,
$$\varphi_{\mu_{i-1}}(\overline{1^{k_i}}h)=\displaystyle \sum_{\eta \in {\cal L}_{\mu_{i-1}}^{k_i}} \varepsilon_{\mu_{i-1} \eta} \phi_{\mu_{i-1} | \eta}(1^{k_i}) \varphi_{\eta}(h).$$
On the other hand, by repeated use of the Murnaghan-Nakayama Rule, we have, for all $h \in S_{\ell}$,
$$\varphi_{\mu_{i-1}}(\overline{1^{k_i}}h)=\displaystyle \sum_{\eta \in {\cal L}_{\mu_{i-1}}^{k_i}} \left( \sum_{P \in {\cal P}_{\eta}^{\mu_{i-1}}} (-1)^{L_P} \right) \varphi_{\eta}(h).$$
Hence, for all $\eta \in {\cal L}_{\mu_{i-1}}^{k_i}$,
$$\varepsilon_{\mu_{i-1} \eta} \phi_{\mu_{i-1} | \eta}(1^{k_i})=\displaystyle  \sum_{P \in {\cal P}_{\eta}^{\mu_{i-1}}} (-1)^{L_P}.$$
Similarly, for all $\eta \in {\cal L}_{\lambda_{i-1}}^{k_i}$,
$$\varepsilon_{\lambda_{i-1} \eta} \phi_{\lambda_{i-1} | \eta}(1^{k_i})=\displaystyle  \sum_{P \in {\cal P}_{\eta}^{\lambda_{i-1}}} (-1)^{L_P}.$$
We therefore have, for any $1 \leq i \leq F$,
$$ {\cal B}_{\lambda_i,\mu_i}^{\lambda_{i-1},\mu_{i-1}}= \displaystyle \sum_{P_{\lambda_i} \in {\cal P}_{\lambda_i}^{\lambda_{i-1}}, \, P_{\mu_i} \in {\cal P}_{\mu_i}^{\mu_{i-1}} } \frac{(-1)^{L_{P_{\lambda_i}}} (-1)^{L_{P_{\mu_i}}} }{k_i! (q^d-1)^{k_i} }.$$

Coming back to our expression for $ <\chi^{\lambda}, \, \chi^{\mu}>_{d-sing}$, we see that, for each $1 \leq k \leq w$, once $k_1, \, \ldots , \, k_F$ are chosen, then all the paths obtained by removing a sequence of $k$ $d$-hooks from $\lambda$ and $\mu$ appear exactly once in the sum. We get that $ {\cal A}_{k_1,\ldots,k_F}^{\lambda,\mu}$ is equal to
$$\displaystyle \frac{1}{k_1! \ldots k_F!(q^d-1)^k} \; \; \sum_{\scriptstyle \lambda' \in {\cal L}_{\lambda}^k \atop \scriptstyle \mu' \in {\cal L}_{\mu}^k}  \; \; \sum_{\scriptstyle P_{\lambda} \in {\cal P}_{\lambda'}^{\lambda} \atop \scriptstyle P_{\mu} \in {\cal P}_{\mu'}^{\mu} } (-1)^{L_{P_{\lambda}}} (-1)^{L_{P_{\mu}}} <\chi^{\lambda'}, \, \chi^{\mu'}>_{ d-reg}.$$
Now, for all $\lambda' \in {\cal L}_{\lambda}^k$ and $\mu' \in {\cal L}_{\mu}^k$, $\lambda'$ and $\mu'$ have $d$-weight $w-k$ and $d$-core $\gamma$, and $\mu'$ is simple and disjoint from $\lambda'$. By the induction hypothesis, we thus have
$$ <\chi^{\lambda'}, \, \chi^{\mu'}>_{G_{n-kd}, \, d-reg}= \displaystyle \frac{ (-1)^{w-k} F^{w-k}}{(w-k)! (q^d-1)^{w-k}} \sum_{P_{\lambda'} \in {\cal P}_{\gamma}^{\lambda'}, \, P_{\mu'} \in {\cal P}_{\gamma}^{\mu'} } (-1)^{L_{P_{\lambda'}}} (-1)^{L_{P_{\mu'}}}.$$
Finally, we get
$$ <\chi^{\lambda}, \, \chi^{\mu}>_{d-sing}=
 \displaystyle \sum_{k=1}^w \; \; \sum_{k_1 + \cdots +k_F=k}  {\cal A}_{k_1,\ldots,k_F}^{\lambda,\mu},$$
where $${\cal A}_{k_1,\ldots,k_F}^{\lambda,\mu}=\displaystyle \frac{ (-1)^{w-k} F^{w-k}}{k_1! \ldots k_F! (w-k)! (q^d-1)^w}\; \;  \sum_{P_{\lambda} \in {\cal P}_{\gamma}^{\lambda}, \, P_{\mu} \in {\cal P}_{\gamma}^{\mu} } (-1)^{L_{P_{\lambda}}} (-1)^{L_{P_{\mu}}}.$$
In order to prove our result, it therefore suffices to show that
$$ \displaystyle \sum_{k=1}^w \; \; \sum_{k_1 + \cdots +k_F=k} \; \; \frac{ (-1)^{w-k} F^{w-k}}{k_1! \ldots k_F! (w-k)!} = \frac{ (-1)^{w+1} F^w}{w!}.$$
We rewrite the second sum in the left hand side as a sum over the
partitions of $k$. When doing this, we ``break'' the ordering on the
elements of ${\cal F}_d$, so that we have to count how many times
each given partition appears once its parts are fixed, and divide by
the number of times the corresponding conjugacy class of
$d$-elements appears. More precisely, to each partition $(k^{r_k},
\, \ldots , \, 1^{r_1})$ of $k$ with $r_1+\cdots+r_k=r$ parts, there
are $F(F-1) \cdots (F-r+1)$ ways to associate one distinct
polynomial of ${\cal F}_d$ to each part (which is why we supposed
that $F \geq w$), but (since the order of the polynomials doesn't
matter, only their multiplicity (i.e. $\left( \begin{array}{cc} f &
\\ & g
\end{array} \right )$ and $\left( \begin{array}{cc} g & \\ & f
\end{array} \right )$ are conjugate)) each conjugacy class of
$d$-element of type  $(k^{r_k}, \, \ldots , \, 1^{r_1})$ will appear
$r_1!\ldots r_k!$ times in this way. The left hand side thus becomes
$$\displaystyle \sum_{k=1}^w \frac{(-1)^{w-k} F^{w-k}}{(w-k)!} \sum_{r=1}^k \; \; \sum_{ \scriptstyle (k^{r_k}, \ldots,1^{r_1}) \vdash k \atop \scriptstyle r_1+ \cdots +r_k=r} \;  \frac{F(F-1) \ldots (F-r+1)}{ (1!)^{r_1} r_1! \ldots (k!)^{r_k} r_k!}.$$
To conclude, we just need the following:

\begin{lem}
For any integers $k \geq 1$ and $F \geq k$, we have
$$ \displaystyle \sum_{r=1}^k  \; \; \sum_{ \scriptstyle (k^{r_k}, \ldots,1^{r_1}) \vdash k \atop \scriptstyle r_1+ \cdots +r_k=r} \frac{F(F-1) \ldots (F-r+1)}{ (1!)^{r_1} r_1! \ldots (k!)^{r_k} r_k!}= \frac{F^k}{k!}.$$
\end{lem}

\begin{proof}
The coefficient of $z^k$ in the power series of $exp(z)^F$ is
$F^k/k!$. Indeed,
$$exp(z)^F=exp(Fz)=\displaystyle \sum_{n \geq 0} \frac{F^n}{n!}z^n.$$
On the other hand, we have
$$\displaystyle \left( \sum_{n \geq 0} \frac{1}{n!}z^n \right)^F= \sum_{r_0+ r_1+ \cdots =F} \frac{C_{(r_0, r_1, \ldots)}}{(0!)^{r_0} (1!)^{r_1} \ldots } z^{\sum_{i \geq 0} ir_i},$$
where $C_{(r_0, r_1, \ldots)}$ is the multinomial coefficient $\frac{F!}{r_0!r_1!\ldots}$. The coefficient of $z^k$ in this expression therefore corresponds to all the uplets $(r_0, \, r_1, \, \ldots)$ such that $r_0+ r_1+ \cdots =F$ and $\sum_{i \geq 0} ir_i=\sum_{i>0}ir_i=k$. Hence $(1^{r_1}, \, \ldots) \vdash k$ is a partition with $\sum_{i>0}r_i=:r$ parts, $r_0=F-\sum_{i>0}r_i=F-r$, and $C_{(r_0, r_1, \ldots)}=\frac{F!}{(F-r)!r_1!\ldots}$. Finally, we get that the coefficient of $z^k$ is
$$ \displaystyle \sum_{r=1}^k \; \; \sum_{\scriptstyle (k^{r_k}, \ldots,1^{r_1}) \vdash k \atop \scriptstyle r_1+ \cdots +r_k=r} \frac{F(F-1) \ldots (F-r+1)}{ (1!)^{r_1} r_1! \ldots (k!)^{r_k} r_k!},$$
which is the desired result.
\end{proof}

\noindent
{\bf{Remark:}} In the proof, we have used the fact that $F$ is an integer. However, the equality given for each integer $k \geq 1$ in the statement of the lemma is an equality between two polynomial functions in $F$, satisfied for an infinite number of values of $F$. It is therefore in fact an equality in the polynomial ring ${\bf Q}[F]$.

\medskip
Using the previous lemma, we get
$$ \begin{array}{cl} \displaystyle \sum_{k=1}^w \; \; \sum_{k_1 + \cdots +k_F=k} \; \frac{ (-1)^{w-k} F^{w-k}}{k_1! \ldots k_F! (w-k)!} & = \displaystyle \sum_{k=1}^w  \frac{ (-1)^{w-k} F^{w-k}}{ (w-k)!} \frac{F^k}{k!} \\ &  =\displaystyle  \frac{F^w}{w!} \sum_{k=1}^w \frac{ (-1)^{w-k} w!}{k! (w-k)!} \\ & = \displaystyle \frac{(-1)^w F^w}{w!} \sum_{k=1}^w (-1)^k \left( \begin{array}{c} w \\k \end{array} \right) \\ & =  \displaystyle \frac{(-1)^w F^w}{w!} [ (1-1)^w -1 ] \\ & = \displaystyle \frac{ (-1)^{w+1} F^w}{w!}. \end{array}$$

Putting everything back together, we finally obtain
$$<\chi^{\lambda}, \, \chi^{\mu}>_{d-sing}= \displaystyle \frac{ (-1)^{w+1} F^w}{w! (q^d-1)^w} \; \sum_{P_{\lambda} \in {\cal P}^{\lambda}_{\gamma}, \, P_{\mu} \in {\cal P}^{\mu}_{\gamma}} \; (-1)^{L_{P_{\lambda}}} (-1)^{L_{P_{\mu}}},$$
which is equivalent to
$$<\chi^{\lambda}, \, \chi^{\mu}>_{d-reg}= \displaystyle \frac{ (-1)^w F^w}{w! (q^d-1)^w} \; \sum_{P_{\lambda} \in {\cal P}^{\lambda}_{\gamma}, \, P_{\mu} \in {\cal P}^{\mu}_{\gamma}} \; (-1)^{L_{P_{\lambda}}} (-1)^{L_{P_{\mu}}}.$$

\medskip

To conclude by induction, it only remains to study the case $w=1$.
This case can in fact be treated without the simplicity hypothesis
on $\mu$.

We take $\lambda$ and $\mu$, distinct partitions of $n$ of
$d$-weight 1, and with the same $d$-core $\gamma$. Then ${\cal
L}_\lambda^1={\cal L}_{\mu}^1=\{ \gamma \}$, and ${\cal
L}_\lambda^i={\cal L}_{\mu}^i=\emptyset$ for $i>1$, and we write
${\cal P}_{\gamma}^{\lambda}=\{ P_{\lambda} \}$ and ${\cal
P}_{\gamma}^{\mu}=\{ P_{\mu} \}$. As above, it suffices to consider
the inner product of $\chi^{\lambda}$ and $\chi^{\mu}$ on the set of
$d$-singular elements of $d$-weight at most 1. The minimal
polynomial of any such element has exactly one irreducible factor of
degree $d$, and it has multiplicity one. Any $d$-singular element $g
\in G$ of weight 1 can thus be written uniquely $g=\rho
\sigma=\sigma \rho$, where $\rho$ can be seen as a $d$-element of
$G_d=GL(d,q)$, and $\sigma$ as a $d$-regular element of
$G_{n-d}=GL(n-d,q)$.

We have
$$\begin{array}{cl} <\chi^{\lambda}, \, \chi^{\mu}>_{d-sing} &  =
\displaystyle \sum_{ \scriptstyle \rho \in G_d/\sim \atop \scriptstyle \rho \; d-sing}
\; \; \sum_{ \scriptstyle \sigma \in G_{n-d}/ \sim \atop \scriptstyle \sigma \; d-reg}
\frac{1}{|C_{G_n}(\rho \sigma)|} \chi^{\lambda}(\rho \sigma)
\chi^{\mu}(\rho \sigma) \\
 & = \displaystyle \sum_{ \scriptstyle \rho \in G_d/\sim \atop \scriptstyle \rho \;  d-sing}
 \frac{a_{\lambda \gamma}^{\rho} a_{\mu \gamma}^{\rho}}{|C_{G_d}(\rho )|}
\sum_{ \scriptstyle \sigma \in G_{n-d}/\sim \atop \scriptstyle \sigma \; d-reg}
\frac{\chi^{\lambda}(\rho \sigma) \chi^{\mu}(\rho \sigma)}{|C_{G_{n-d}}(\sigma)|} \\
 &  =
\displaystyle  F  \sum_{\nu \vdash 1} \frac{a_{\lambda
\gamma}^{\rho^{\nu}} a_{\mu
\gamma}^{\rho^{\nu}}}{|C_{G_d}(\rho^{\nu} )|} <\chi^{\gamma}, \,
\chi^{\gamma}>_{G_{n-d},d-reg}.
 \end{array}$$
By Lemma 4.7, we thus have
$$ \begin{array}{cl} <\chi^{\lambda}, \, \chi^{\mu}>_{d-sing} &  = \displaystyle F \frac{\varepsilon_{\lambda  \gamma} \varepsilon_{\mu \gamma} \phi_{\lambda | \gamma}(1) \phi_{\mu | \gamma}(1)}{ |C_{S_1}(1)| \; |T_{\overline{(1)}}|}  <\chi^{\gamma}, \, \chi^{\gamma}>_{G_{n-d},d-reg} \\  &  = \displaystyle F \frac{\varepsilon_{\lambda  \gamma} \varepsilon_{\mu \gamma} \phi_{\lambda | \gamma}(1) \phi_{\mu | \gamma}(1)}{ q^d-1}  <\chi^{\gamma}, \, \chi^{\gamma}>_{G_{n-d},d-reg} \\ & = \displaystyle F \frac{\varepsilon_{\lambda  \gamma} \varepsilon_{\mu \gamma} \phi_{\lambda | \gamma}(1) \phi_{\mu | \gamma}(1)}{ q^d-1}, \end{array}$$
since $\gamma$ is a $d$-core, whence vanishes on $d$-singular
elements of $GL(n-d,q)$. Finally, as after Lemma 4.8, we have
$\varepsilon_{\lambda  \gamma} \phi_{\lambda | \gamma}(1)=
(-1)^{L_{P_{\lambda}}}$ and $ \varepsilon_{\mu \gamma} \phi_{\mu |
\gamma}(1)= (-1)^{L_{P_{\mu}}}$, which gives us
$$ <\chi^{\lambda}, \, \chi^{\mu}>_{d-sing}  = \displaystyle \frac{F}{q^d-1}(-1)^{L_{P_{\lambda}}} (-1)^{L_{P_{\mu}}}.$$

Hence the result is true for $w=1$. Thus, by induction, we have
that, for any $w \geq 1$, if $\lambda$ and $\mu$ are partitions of
$n$ of $d$-weight $w$, with the same $d$-core $\gamma$, and if $\mu$
is simple and disjoint from $\lambda$, then
$$<\chi^{\lambda}, \, \chi^{\mu}>_{d-reg}= \displaystyle \frac{ (-1)^w F^w}{w! (q^d-1)^w} \sum_{P_{\lambda} \in {\cal P}^{\lambda}_{\gamma}, \, P_{\mu} \in {\cal P}^{\mu}_{\gamma}} (-1)^{L_{P_{\lambda}}} (-1)^{L_{P_{\mu}}}.$$

To conclude, we just note that, in fact, $(-1)^{L_{P_{\lambda}}}$
(resp. $(-1)^{L_{P_{\mu}}}$) is independent on the path $P_{\lambda}
\in {\cal P}^{\lambda}_{\gamma}$ (resp. $P_{\mu} \in {\cal
P}^{\mu}_{\gamma}$) (cf \cite{James-Kerber}, Theorem 2.7.27). We
therefore write $(-1)^{L_{P_{\lambda}}}=\varepsilon_{\lambda}$ for
any $P_{\lambda} \in {\cal P}^{\lambda}_{\gamma}$ (resp.
$(-1)^{L_{P_{\mu}}}=\varepsilon_{\mu}$, for any $P_{\mu} \in {\cal
P}^{\mu}_{\gamma}$). We get
$$<\chi^{\lambda}, \, \chi^{\mu}>_{d-reg}= \displaystyle \frac{ (-1)^w F^w}{w! (q^d-1)^w} |{\cal P}^{\lambda}_{\gamma}| \, |{\cal P}^{\mu}_{\gamma}| \varepsilon_{\lambda} \varepsilon_{\mu} \neq 0.$$
This ends the proof of Theorem 4.6.

\subsection{Nakayama Conjecture for blocks of small weight}

We are now in position to prove the second direction of the Nakayama
Conjecture for unipotent $d$-blocks of $G$, provided their weight is
``small'' compared to $d$. For partitions $\lambda$ and $\mu$ of
$n$, we will write $\lambda \sim \mu$ if $\chi_{\lambda}$ and
$\chi_{\mu}$ are directly linked across $d$-singular elements, and
$\lambda \equiv \mu$ if $\chi_{\lambda}$ and $\chi_{\mu}$ belong to
the same unipotent $d$-block of $G$. First, take $\lambda$ and $\mu$
of weight $w \neq 0$, with the same $d$-core, and suppose that $3w
\leq d$ and $F \geq w$. Then there exists $\nu \vdash n$ which is
simple and disjoint from both $\lambda$ and $\mu$. Thus, by Theorem
4.6, $\lambda \sim \nu$ and $\mu \sim \nu$, so that $\lambda \equiv
\mu$. In fact, we can do better than that:

\begin{thm}
Let $d>0$ be an integer. Suppose $w >2$ is an integer, and $w \leq
F$ (the number of irreducible monic polynomials of degree $d$ over
${\bf F}_q$). If $d \geq 2w-1$, then the unipotent $d$-blocks of $G$
of weight $w$ satisfy the Nakayama Conjecture.
\end{thm}

\begin{proof}
Take $\lambda$ and $\mu$ partitions of $n$ of weight $w$, with the
same $d$-core $\gamma$. The proof is easy to understand on a picture
of the abacus, even though it is difficult to write down.

Throughout the proof, we will use the following notations: we say
that $\lambda$ {\emph{uses}} a runner of the abacus if some $d$-hook
removal of $\lambda$ is stored on this runner; we say that $\mu
\subset \lambda$ if all the runners used by $\mu$ are used by
$\lambda$.

First, suppose that $\lambda$ and $\mu$ use each at most $w-1$ runners. Then there exists $\nu \vdash n$ which is simple and disjoint from $\lambda$, and of $d$-core $\gamma$. By Theorem 4.6, we have $\lambda \sim \nu$. Now
\begin{itemize}
\item{}
if $\mu \not \subset \nu$, then there exists $\zeta \vdash n$,
disjoint from $\nu$, with $d$-core $\gamma$, and which uses a single
runner, which runner is also used by $\mu$. Then $\nu \sim \zeta$.
Now there exists $\xi$, which is simple, disjoint from $\mu$ (and
thus from $\zeta$), and with $d$-core $\gamma$, so that $\zeta \sim
\xi$ and $\xi \sim \mu$. Hence $\lambda \equiv \mu$.
\item{}
if $\mu \subset \nu$, then there exists $\zeta$, disjoint from
$\nu$, with $d$-core $\gamma$, and which uses a single runner, which
runner is not used by $\mu$. Thus $\nu \sim \zeta$. There exists
$\xi$, simple, disjoint from $\zeta$, with $d$-core $\gamma$, and
such that at least one runner used by $\mu$ is not used by $\xi$; we
have $\zeta \sim \xi$. Then $\xi \sim \delta$, where $\delta$ has
$d$-core $\gamma$, and uses a single runner, used by $\mu$ but not
by $\xi$. Now $\delta \sim \eta$, simple, with $d$-core $\gamma$,
and disjoint from $\mu$ (and thus from $\delta$). Thus $\eta \sim
\mu$, and $\lambda \equiv \mu$.
\end{itemize}

Next, suppose that $\lambda$ is simple (i.e. lies on $w$ runners). If $\mu$ is disjoint from $\lambda$, then, by Theorem 4.6, $\lambda \sim \mu$, so that $\lambda \equiv \mu$. Suppose thus that $\lambda$ and $\mu$ share a runner.
\begin{itemize}
\item{}
Suppose $\mu$ lies on at most $w-1$ runners.\\
If $\mu \subset \lambda$, then $\lambda \sim \nu$ with $d$-core
$\gamma$, disjoint from $\lambda$, and using a single runner. Now
$\nu \sim \zeta$, simple, with $d$-core $\gamma$, and such that $\mu
\not \subset \zeta$. Then $\zeta \sim \xi$ with $d$-core $\gamma$,
using a single runner, which runner is used by $\mu$. Then $\xi \sim
\eta$ simple, with $d$-core $\gamma$, and disjoint from $\mu$ (and
thus from $\xi$), so that $\eta \sim \mu$. Hence $\lambda \equiv
\mu$.

If $\mu \not \subset \lambda$, then $\lambda \sim \nu$ with $d$-core
$\gamma$, disjoint from $\lambda$, and using a single runner, which
runner is used by $\mu$. Thus $\nu \sim \zeta$ simple, with $d$-core
$\gamma$, and disjoint from $\mu$ (and thus from $\nu$), so that
$\zeta \sim \mu$. Hence $\lambda \equiv \mu$.

\item{}
Suppose $\mu$ is simple.\\
If there exists a runner which is used by neither $\lambda$ nor
$\mu$, then there exists $\nu$ with $d$-core $\gamma$ which lies on
this single runner (and thus is disjoint from both $\lambda$ and
$\mu$), so that $\lambda \sim \nu$ and $\nu \sim \mu$, whence
$\lambda \equiv \mu$.

If not, then, necessarily, $d=2w-1$, and $\lambda$ and $\mu$ share a
single runner. Since $w>2$, there is at least one runner used by
$\mu$ and not by $\lambda$. Then $\lambda \sim \nu$ with $d$-core
$\gamma$ which lies on this single runner (and is thus disjoint from
$\lambda$). Now $\nu \sim \zeta$ simple, with $d$-core $\gamma$,
disjoint from $\nu$ and which doesn't use at least one runner which
is used by $\lambda$ but not by $\mu$. Then $\zeta \sim \xi$ with
$d$-core $\gamma$ and which lies on this single runner. In
particular, $\xi$ is disjoint from $\mu$, so that $\xi \sim \mu$
(since $\mu$ is simple). Hence $\lambda \equiv \mu$.
\end{itemize}
This concludes the proof.
\end{proof}

\medskip
\noindent
{\bf{Remarks:}}

\noindent 1. The case of blocks of weight 0 is easily dealt with. A
partition of $d$-weight 0 is its own $d$-core, so that the
corresponding unipotent character is alone in its unipotent
$d$-block. Hence, for any $d>0$, the unipotent $d$-blocks of $G$ of
weight 0 satisfy the Nakayama Conjecture.

\noindent 2. As we have remarked earlier, in the case $w=1$, we
don't need to use simple partitions. For any $d >0$ and for any two
partitions $\lambda$ and $\mu$ of $n$ of $d$-weight 1 and with the
same $d$-core, the unipotent characters $\chi_{\lambda}$ and
$\chi_{\mu}$ are directly linked across $d$-singular elements. Hence
the unipotent $d$-blocks of weight 1 of $G$ satisfy the Nakayama
Conjecture.

\noindent 3. As can be seen in the proof of Theorem 4.10, the case
$w=2$ is problematic only in the case $d=2w-1$, and the problem can
easily be solved if we have another runner we can use. Hence the
result of the theorem is also true for $w=2$ if we only suppose $2w
\leq d$.

Furthermore, the case $w=2$ (which is fairly small) can in fact be
studied on its own, and one can show that, for any $d>0$, the
unipotent $d$-blocks of $G$ of weight 2 satisfy the Nakayama
Conjecture (provided $F \geq 2$).

\noindent 4. Since the $d$-weight of a unipotent $d$-block of $G$ is
at most $n/d$, we see that, if $n \leq d(d+1)/2$, then all the
unipotent $d$-blocks of $G$ satisfy the Nakayama Conjecture
(provided $F \geq n/d$).

\subsection{And then?}
%From the two cases we just studied, we see that the first possible
%counter-example to the unipotent $d$-blocks of $GL(n,q)$ satisfying
%an analogue of the Nakayama Conjecture is given by the unipotent
%2-blocks of $GL(4,q)$. However, using the tables given by L. E.
%Dickson (cf \cite{Dickson}, Chapter X, section 223) for the
%conjugacy classes of $GL(4,q)$, as well as the work of R. Steinberg
%(cf \cite{Steinberg}), where the values of the unipotent characters
%of $GL(4,q)$ are given explicitly, we can easily test directly the
%result. It turns out that, here again, any two (irreducible)
%unipotent characters labeled by partitions with the same 2-core are
%directly linked across 2-singular elements. In particular, the
%unipotent 2-blocks of $GL(4,q)$ satisfy an analogue of the Nakayama
%Conjecture.

%\medskip
The proof of the Nakayama Conjecture in the general case seems to be hard. In \cite{KOR}, the authors prove that, if
$B$ is an $\ell$-block of weight $w$ of some symmetric group, then
there is a {\emph{generalized perfect isometry}} between $B$ and the
set of irreducible (complex) characters of the wreath product ${\bf
Z}_{\ell} \wr S_w$. This means that the restricted scalar products
of characters of $B$ across $\ell$-regular elements are, up to a
sign, the same as the restricted scalar products of characters of
${\bf Z}_{\ell} \wr S_w$ across some (carefully chosen) set of
elements (called {\emph{regular}} elements).

It turns out that the computations are easier to carry out in the
wreath product, where it can be shown that every character is
directly linked across regular elements to the trivial character.
Coming back to the symmetric group, one gets that any two characters
in the combinatorial $\ell$-block are linked across $\ell$-regular
elements. This proves that the combinatorial $\ell$-blocks and the
$\ell$-blocks are the same, or, equivalently, that the $\ell$-blocks
satisfy the Nakayama Conjecture.

\smallskip
Now, we might try and do something similar in the case of $GL(V)$.
This would involve finding the relevant wreath product, and the set
of conjugacy classes we want to distinguish in it. Given a
combinatorial unipotent $d$-block of weight $w$, a {\emph{natural}}
candidate for the wreath product is $GL(d,q) \wr S_w$. However, even
in easy cases (like $d=1$), this doesn't seem to contain enough
information (in particular, it is hard to have powers of $q$
appearing).

So far, we can only conjecture that, for any $d>0$, the unipotent
$d$-blocks of $G$ satisfy an analogue of the Nakayama Conjecture.

\section{Second Main Theorem}

One of the very striking properties of {\emph{ordinary}} blocks of
finite groups is Brauer's Second Main Theorem (cf e.g.
\cite{Navarro}). In \cite{KOR}, the authors give an analogue of this
in the context of generalized blocks. Depending on the generalized
sections we define, the blocks we obtain may satisfy this analogue
or not.

\subsection{Domination, Second Main Theorem Property}

For all the definitions and results in this section, we refer to
\cite{KOR}. We start by defining blocks for the centralizers of
$d$-elements. Recall that the $d$-blocks of $G$ are defined by
orthogonality across the set ${\cal Y}_d(1)$ of $d$-regular
elements. As a consequence, they separate ${\cal Y}_d(1)$ from its
complement (cf \cite{KOR}, Corollary 1.2). For $x \in {\cal X}_d$,
we define the $d${\emph{-blocks}} of $C_G(x)$ to be the smallest
(non-empty) subsets of Irr$(C_G(x))$ such that irreducible
characters in distinct subsets are orthogonal across $x{\cal
Y}_d(x)$. We can equally define them to be non-empty subsets of
Irr$(C_G(x))$ which are minimal subject to separating ${\cal
Y}_d(x)$ from its complement in $C_G(x)$ (or, equivalently, to
separating $x{\cal Y}_d(x)$ from its complement (since $x$ is
central in $C_G(x)$)).

\medskip
We now turn to the notion of domination. Suppose $\chi \in
\mbox{Irr}(G)$ and $\beta$ is a union of $d$-blocks of $C_G(x)$ for
some $x \in {\cal X}_d$. We define a generalized character
$\chi^{(\beta)}$ of $C_G(x)$ via
$$\chi^{(\beta)}= \displaystyle \sum_{\mu \in \beta} < \mbox{Res}^G_{C_G(x)}(\chi),\mu> \mu.$$

\begin{defn}
Let $x \in {\cal X}_d$ and $b$ be a $d$-block of $C_G(x)$. We say
that a $d$-block $B$ of $G$ {\bf{dominates}} $b$ if there exist
$\chi \in B$ and $y \in {\cal Y}_d(x)$ such that $\chi^{(b)}(xy)
\neq 0$.
\end{defn}

We see that, for $x \in {\cal X}_d$, if $\chi \in B$ for some
$d$-block $B$ of $G$, then, for each $y \in {\cal Y}_d(x)$, we have
$\chi(xy)=\sum_b \chi^{(b)}(xy)$, where $b$ runs through the set of
$d$-blocks of $C_G(x)$ dominated by $B$.

\medskip
Note that, for $x \in {\cal X}_d$, each $d$-block of $C_G(x)$ is
dominated by at least one $d$-block of $G$.

\begin{defn}
We say that the $d$-blocks of $G$ satisfy the Second Main Theorem
Property if, for each $x \in {\cal X}_d$ and each $d$-block $b$ of
$C_G(x)$, $b$ is dominated by a unique $d$-block of $G$.
\end{defn}

Note that, if, instead of ${\cal X}_d$, we take the set of
$r$-elements of $G$ ($r$ a prime), and, instead of ${\cal Y}_d(x)$,
we take the set of $r$-regular elements of $C_G(x)$, then we obtain
the $r$-blocks of $G$, and they do satisfy the Second Main Theorem
Property.

\medskip
Using the fact that, for $x \in {\cal X}_d$, irreducible characters
of $C_G(x)$ in distinct $d$-blocks are orthogonal across $x{\cal
Y}_d(x)$, one proves easily the following:

\begin{prop}(\cite{KOR}, Corollary 2.2)
The $d$-blocks of $G$ satisfy the Second Main Theorem Property if
and only if, for each $d$-block $B$ of $G$, there is, for each $x
\in {\cal X}_d$, a (possibly empty) union $\beta(x,B)$ of $d$-blocks
of $C_G(x)$ such that, for each irreducible character $\chi \in B$
and each character $\mu \in \beta(x,B)$, we may find a complex
number $c_{\chi, \mu}$ such that, for each $y \in {\cal Y}_d(x)$, we
have
$$\chi(xy)=\displaystyle \sum_{\mu \in \beta(x,B)} c_{\chi, \mu} \mu(xy),$$
and, furthermore, $\beta(x,B)$ and $\beta(x,B')$ are disjoint
whenever $B$ and $B'$ are distinct $d$-blocks of $G$.
\end{prop}

The following theorem enlighten the link between the Second Main
Theorem Property and Brauer's Second Main Theorem:
\begin{thm} (\cite{KOR}, Corollary 2.3)
Suppose that the $d$-blocks of $G$ satisfy the Second Main Theorem
property. Then:

\noindent (i) Irreducible characters of $G$ which are in distinct
$d$-blocks are orthogonal across {\bf{each}} ${\cal Y}_d$-section of
$G$.

\noindent (ii) If $x \in {\cal X}_d$ and $\sum_{\chi \in
\mbox{Irr}(G)}a_{\chi} \chi$ is a class function which vanishes
identically on the ${\cal Y}_d$-section of $x$ in $G$, then, for
each $d$-block $B$ of $G$, $\sum_{\chi \in B} a_{\chi} \chi$ also
vanishes identically on the ${\cal Y}_d$-section of $x$ in $G$.

\noindent (iii) $d$-blocks of $G$ separate ${\cal Y}_d$-sections of
$G$.
\end{thm}

\subsection{Second Main Theorem Property for combinatorial unipotent
$d$-blocks}

We first study the $d$-blocks of the centralizers of the
$d$-elements. If $x \in {\cal X}_d$ has $d$-type ${\bf km}$, then,
writing $l=n-{\bf km}d$, we have $C_G(x)=H_0 \times H_1$ where $H_1
\cong GL(l,q)$ and $H_0 \leq G_0 \cong GL({\bf km}d,q)$. Then
Irr$(C_G(x))=$Irr$(H_0) \otimes$Irr$(H_1)$. Note that, as we noted
before, we may consider $x$ as an element of $H_0$, and then $x{\cal
Y}_d(x)=\{ (x,y) \in H_0 \times H_1, \, y \in {\cal Y}_d^l(1) \}$.

Take $\chi_0, \psi_0 \in \mbox{Irr}(H_0)$ and $\chi_1, \psi_1 \in
\mbox{Irr}(H_1)$. We have
$$\begin{array}{cl} <\chi_0 \otimes \chi_1, \psi_0 \otimes \psi_1>_{x{\cal Y}_d(x)} &
= \displaystyle \frac{1}{|C_G(x)|} \sum_{y \in {\cal Y}_d(x)}
(\chi_0 \otimes \chi_1)(xy) \overline{( \psi_0 \otimes \psi_1)(xy)}
\\ & =\displaystyle \frac{1}{|C_G(x)|} \sum_{y \in {\cal Y}_d^l(1)}
\chi_0(x)  \chi_1(y) \overline{ \psi_0(x)} \overline{ \psi_1(y)} \\
&=\displaystyle \frac{\chi_0(x) \overline{\psi_0(x)}}{|C_G(x)|}
\sum_{y \in {\cal Y}_d^l(1)}   \chi_1(y)  \overline{ \psi_1(y)}  \\
&=\displaystyle \frac{\chi_0(x) \overline{\psi_0(x)}}{|H_0|}
<\chi_1,\psi_1>_{{\cal Y}_d^l(1)}. \end{array}$$ Since $x$ is
central in $H_0$, we have $\chi_0(x) \overline{\psi_0(x)} \neq 0$,
and we see that $\chi_0 \otimes \chi_1$ and $ \psi_0 \otimes \psi_1$
are directly $x{\cal Y}_d(x)$-linked if and only if $\chi_1$ and
$\psi_1$ are directly ${\cal Y}_d^l(1)$-linked. Extending by
transitivity, we obtain that the $d$-blocks of $C_G(x)$ are the
Irr$(H_0) \otimes b_i$'s, where $b_i$ runs through the set of ${\cal
Y}_d^l(1)$-blocks (i.e. $d$-blocks) of $H_1 \cong GL(l,q)$.

\medskip
In analogy with this, we define the unipotent $d$-blocks and
combinatorial unipotent $d$-blocks of $C_G(x)$ to be the Irr$(H_0)
\otimes b_i$'s, where $b_i$ runs through the sets of unipotent
$d$-blocks and combinatorial unipotent $d$-blocks of $H_1 \cong
GL(l,q)$ respectively.

\medskip
We can now prove that the combinatorial unipotent $d$-blocks of $G$
satisfy the Second Main Theorem Property.

Take any $x \in {\cal X}_d$, and write $C_G(x)=H_0 \times H_1$ as
above. For any combinatorial unipotent $d$-block $B$ of $G$, labeled
by the $d$-core $\gamma$, we set $\beta(x,B)=\mbox{Irr}(H_0) \otimes
b$, where $b$ is the combinatorial unipotent $d$-block of $H_1$
labeled by $\gamma$. For any $\chi_{\mu} \in B$ and $\psi_0 \otimes
\psi_{\lambda} \in \mbox{Irr}(H_0) \otimes b$, we set
$$c_{\chi_{\mu},\psi_0 \otimes \psi_{\lambda}}= \left\{ \begin{array}{ll} \alpha^{x}_{\mu \lambda}
 & \mbox{if} \; \psi_0=1_{H_0} \\ 0 & \mbox{otherwise} \end{array} \right. ,$$
where the $ \alpha_{\mu \lambda}^{x}$'s are the MN-coefficients,
obtained from the Murnaghan-Nakayama rule for unipotent characters.

Then the definition of the $ \alpha_{\mu \lambda}^{x}$'s shows that,
for each $x \in {\cal X}_d$, the $\beta(x,B)$'s and
$c_{\chi_{\mu},\psi_0 \otimes \psi_{\lambda}}$'s satisfy the
hypotheses of Proposition 5.3. Indeed, for each combinatorial
unipotent $d$-block $B$ of $G$, for each $\chi_{\mu} \in B$ and for
each $y \in {\cal Y}_d(x)$, we have
$$\chi_{\mu}(xy)=\displaystyle \sum_{\psi_0 \otimes \psi_{\lambda} \in \beta(x,B)}
c_{\chi_{\mu},\psi_0 \otimes \psi_{\lambda}} (\psi_0 \otimes
\psi_{\lambda})(xy),$$ and, furthermore, $\beta(x,B)$ and
$\beta(x,B')$ are disjoint whenever $B$ and $B'$ are distinct
combinatorial unipotent $d$-blocks of $G$. This implies that, for
each $x \in {\cal X}_d$, each combinatorial unipotent $d$-block of
$C_G(x)$ is dominated by a unique combinatorial unipotent $d$-block
of $G$. This proves the following:

\begin{thm}
For any positive integers $n$ and $d$, the combinatorial unipotent
$d$-blocks of $G$ satisfy the Second Main Theorem Property. In
particular, if $d=1$, then the unipotent $d$-blocks of $G$ satisfy
the Second Main Theorem Property, and, for any $d>0$, the unipotent
$d$-blocks of $G$ of weight 0, 1, 2, and at most $(d+1)/2$ satisfy
the Second Main Theorem Property.
\end{thm}

\bigskip
\noindent {\bf{Remark:}} It is easy to see that one can change the
definition of $d$-regular element to ``an element whose minimal
polynomial has no irreducible factor of degree $d$ (except maybe
$X-1$)'' (instead of ``degree divisible by $d$'') without affecting
{\bf{any}} of the results we have proved, apart from Theorem 4.5
(Nakayama Conjecture in the case $d=1$).

\bigskip
\noindent {\bf{Acknowledgements:}}

\medskip
I wish to thank deeply Meinolf Geck, Radha Kessar, Attila Mar\'oti
and Geoffrey R. Robinson for the many discussions we had about this
work. I also wish to thank Paul Fong and Bhama Srinivasan for their
kindness and support. Finally, I thank Olivier Mathieu for his help
with the proof of Lemma 4.9.

\end{document}